\documentclass[a4paper,reqno,10pt,oneside]{amsart}
\usepackage{amsmath}
\usepackage[left=2.61cm,right=2.61cm,top=2.72cm,bottom=2.72cm]{geometry}
\usepackage[colorlinks=true,urlcolor=blue,
citecolor=red,linkcolor=blue,linktocpage,pdfpagelabels,
bookmarksnumbered,bookmarksopen]{hyperref}
\usepackage[english]{babel}
\usepackage{graphicx}
\usepackage[toc,page,header]{appendix}
\usepackage{mathrsfs}
\usepackage{amssymb,amsmath, amsfonts, amsthm}
\usepackage{latexsym}
\usepackage{enumitem}

\allowdisplaybreaks[1]

\numberwithin{equation}{section}

\DeclareMathOperator{\divergence}{div}

\renewcommand{\(}{\left(}
\renewcommand{\)}{\right)}
\renewcommand{\[}{\left[}
\renewcommand{\]}{\right]}

\newtheorem{theorem}{Theorem}[section]
\newtheorem{proposition}[theorem]{Proposition}
\newtheorem{corollary}[theorem]{Corollary}
\newtheorem{lemma}[theorem]{Lemma}
\theoremstyle{definition}

\theoremstyle{definition}

\theoremstyle{definition}

\renewcommand{\le}{\leqslant}
\renewcommand{\ge}{\geqslant}

\newcommand{\Cl}{{\mathcal L}}

\renewcommand{\S}{\mathbb{S}}

\newcommand{\beq}{\begin{equation}}
\newcommand{\eeq}{\end{equation}}
\newcommand{\beqs}{\begin{equation*}}
\newcommand{\eeqs}{\end{equation*}}
\newcommand{\beqn}{\begin{eqnarray}}
\newcommand{\eeqn}{\end{eqnarray}}
\newcommand{\beqns}{\begin{eqnarray*}}
\newcommand{\eeqns}{\end{eqnarray*}}
\newcommand{\bdoc}{\begin{document}}
\newcommand{\edoc}{\end{document}}
\newcommand{\be}{\begin{enumerate}}
\newcommand{\ee}{\end{enumerate}}
\newcommand{\bdescr}{\begin{description}}
\newcommand{\edescr}{\end{description}}
\newcommand{\ba}{\begin{array}}
\newcommand{\ea}{\end{array}}
\newcommand{\intR}{\int_{\mathbb R^N}}
\newcommand{\R}{\mathbb R}
\newcommand{\RN}{\mathbb{R}^N}
\newcommand{\B}{\mathbb B}
\newcommand{\C}{\mathbb C}
\renewcommand{\H}{\mathcal H}
\renewcommand{\L}{\mathbb L}
\newcommand{\parallelsum}{\mathbin{\!/\mkern-5mu/\!}}
\newcommand{\e}{\varepsilon}
\newcommand{\SD}{\Sigma_D}
 \renewcommand{\(}{\left(}
\renewcommand{\)}{\right)}
\renewcommand{\[}{\left[}
\renewcommand{\]}{\right]}
\renewcommand{\appendixpagename}{\centering Appendix}

\newcommand{\todo}[1]{\text{\colorbox{yellow}{#1}}}

\begin{document}
\title[Classification of extremals of Caffarelli-Kohn-Nirenberg inequalities]{On the classification of extremals of Caffarelli-Kohn-Nirenberg inequalities}

\author{Giulio Ciraolo}
\address{Giulio Ciraolo. Dipartimento di Matematica "Federigo Enriques",
Universit\`a degli Studi di Milano, Via Cesare Saldini 50, 20133 Milano, Italy}
\email{giulio.ciraolo@unimi.it}

\author{Camilla Chiara Polvara}
\address{Camilla Chiara Polvara. Dipartimento di Matematica, Sapienza Universit\`a di Roma,  P.le Aldo Moro 2, 00185 Roma, Italy}
\email{camilla.polvara@uniroma1.it}

%\thanks{\emph{Acknowledgements.} Research partially supported by Gruppo Nazionale per l'Analisi Matematica, la Pro\-ba\-bi\-li\-t\`a e le loro Applicazioni (GNAMPA) of the Istituto Nazionale di Alta Matematica (INdAM)}

\begin{abstract}
We consider a family of critical elliptic equations which arise as the Euler-Lagrange equation of Caffarelli-Kohn-Nirenberg inequalities, possibly in convex cones in $\mathbb{R}^d$, with $d\geq 2$. We classify positive solutions without assuming that the solution has finite energy and when the intrinsic dimension $n \in (\frac{3}{2},5]$.  
\end{abstract}

\maketitle

\section{Introduction}
In this paper we are interested in the classification of extremals of Caffarelli-Kohn-Nirenberg (CKN) type inequalities 
\begin{equation}\label{CKN}
\bigg(\int_{\R^d}\frac{|u|^p}{|x|^{bp}}dx\bigg)^\frac{2}{p}\le C_{a,b}\int_{\R^d}\frac{|\nabla u|^2}{|x|^{2a}}dx \,.
\end{equation}
CKN inequalities \eqref{CKN} were established in \cite{CKN} (see also \cite{ILIN, Lin, Mazya}) and they hold for any $u \in \mathcal{D}^{a,b}(\R^d)$, where
\begin{equation*}
\mathcal{D}^{a,b} (\R^d) : = \left\{ u \in L^p(\R^d,|x|^{-b}dx) \textmd{ such that  } |x|^{-a} |Du| \in L^2(\R^d,dx) \right\} \,,
\end{equation*} 
with $a \leq b\leq 1+a$ for $d \geq 3$ and $a < b\leq 1+a$ for $d = 2$. 
The constant $C_{a,b}$ appearing in \eqref{CKN} is the optimal constant, it depends only on $a,b$ and $d$. The exponent  
$$
p=\frac{2d}{d-2 + 2(b-a)} 
$$
is chosen in order to have the invariance of the inequality under scalings. 
We shall assume that $a<a_c$, where 
$$
a_c:=\frac{d-2}{2} \,,
$$ 
since the case $a>a_c$ is dual to it (see \cite[Theorem 1.4 (ii)]{CW}) and for $a=a_c$ the inequality \eqref{CKN} fails to be true \cite{CKN}. 

We notice that if $a=b=0$ then \eqref{CKN} reduces to the classical Sobolev inequality, and in this case the study of symmetry properties of extremals for \eqref{CKN} goes back to Aubin \cite{A} and Talenti \cite{T}. Classification results for positive solutions to the corresponding Euler-Lagrange equation
\begin{equation} \label{eqLapl}
\Delta u + u^{\frac{d+2}{d-2}} = 0 \quad \textmd{ in } \R^d 
\end{equation}
were also obtained by Gidas, Ni and Nirenberg in \cite{GNN} under the assumption that the solution behaves at infinity as the fundamental solution. The classification result asserts that positive solutions to \eqref{eqLapl} are radially symmetric with respect to some centre, up to rescalings.

When $(a,b) \neq (0,0)$, a relevant property of \eqref{CKN} is that minimizers are not necessarily radially symmetric. In \cite{FS} Felli and Schneider proved that symmetry breaking occurs for some values of $a$ and $b$. More precisely, by setting
\begin{equation} \label{alpha_def}
\alpha = \frac{(1+a-b)(a_c-a)}{a_c-a+b} \,,
\end{equation}
with
\begin{equation} \label{n_def}
n= \frac{d}{1+a-b} \,,
\end{equation}
it was shown in \cite{FS} that minimizers of \eqref{CKN} are not radially symmetric whenever
\begin{equation} \label{alphaFS}
\alpha > \sqrt{\frac{n-1}{d-1}} \,.
\end{equation}
In \cite{DEL}, Dolbeault, Esteban and Loss proved that the range \eqref{alphaFS} is sharp. 
Indeed, the main result of \cite{DEL} is the following: let $u \in \mathcal{D}^{a,b}(\R^d)$ be a positive solution to   
\begin{equation}\label{eq_CKN}
\divergence \left( |x|^{-2a} Du \right) + |x|^{-bp} u^{p-1} = 0 \quad \textmd{ in } \R^d  \,,
\end{equation}
and assume that 
\begin{equation} \label{alpha_optimal}
\alpha \leq \sqrt{\frac{n-1}{d-1}} \,,
\end{equation}
then $u$ is radially symmetric and it is of the form
\begin{equation} \label{u_radial}
u(x) = \left( 1+ |x|^{(p-2)(a_c-a)}\right)^{-\frac{2}{p-2}}  
\end{equation}
for any $x \in \R^d$, up to rescalings. Since \eqref{eq_CKN} is the Euler-Lagrange equation associated to \eqref{CKN} then this result shows that minimizers of \eqref{CKN} are radially symmetric whenever $\alpha$ satisfies \eqref{alpha_optimal} (we will refer to solutions of \eqref{eq_CKN} as critical points of \eqref{CKN}).
We also mention \cite{ABCMP,BBMP,CW,C,L,DE1,DE2,DE3,DE4,DEL2,DELT,DELT2,LW,LW2,SW} where some partial and related results on symmetry of extremals of \eqref{CKN} were obtained before the sharp characterization in \cite{DEL}.

We emphasize that the assumption $u \in \mathcal{D}^{a,b}(\R^d)$ is crucial in \cite{DEL} for obtaining asymptotic estimates on the solution at infinity and at the origin (see also \cite{SV} and \cite{CC} for more general CKN type inequalities). However, this assumption is not necessary when $a=b=0$, as it was shown by Caffarelli, Gidas and Spruck in \cite{CGS} and Li and Zhang in \cite{LZ}, where the classification result was obtained for positive solutions to \eqref{eqLapl} without assuming neither $u\in \mathcal{D}^{0,0}(\R^d)$ nor asymptotic estimates at infinity. Hence, a natural question is whether or not the assumption $u \in \mathcal{D}^{a,b}(\R^d)$ can be removed for the more general equation \eqref{eq_CKN}. This is the main goal of this paper.

As already mentioned, when $a=b=0$ the energy assumption on the solution was removed in \cite{CGS} and \cite{LZ} by using the Kelvin transformation and the methods of moving planes. However, the method of moving planes seems to be helpful for the classification of solutions to \eqref{eq_CKN} only when $a \geq 0$, see Chou and Chu \cite{C} where symmetry results are obtained under the assumption $u \in \mathcal{D}^{a,b}(\R^d)$. Trying to use this approach when $a<0$ and without the energy assumption seems to be a challenging problem. 

In this paper, we prefer to adopt a different approach which is based on integral identities which, in particular, exploits the use of a suitable $P$-function. This approach was recently used by the authors and Farina in \cite{CFP} for classification and rigidity results of \eqref{eqLapl} in complete Riemmanian manifolds with nonnegative Ricci curvature (see also the related papers by Catino and Monticelli \cite{CM} and Fogagnolo, Malchiodi and Mazzieri \cite{FMM}). We also mention that removing the energy assumption for positive solutions to the critical $p$-Laplace equation 
$$
\Delta_p u + u^{p^*-1}= 0 \quad \textmd{ in } \R^d
$$
is a challenging open problem, and some partial results have been recently obtained by Catino, Monticelli and Roncoroni in \cite{CMR}, Ou in \cite{OU} and Vetois in \cite{Vetois}, which have been other sources of inspiration for us.

In this paper, we prove the following result, where we obtain a rigidity result in low \emph{intrinsic} dimension $n$ and without energy assumptions on the solution.

\begin{theorem} \label{thm_main}
Let $u \in \mathcal{D}^{a,b}_{loc} (\R^d)$ be a positive solution to \eqref{eq_CKN} and let $\alpha$ and $n$ be given by \eqref{alpha_def} and \eqref{n_def}, respectively. If 
\begin{equation} \label{alpha_reduced} 
\frac52 < n \leq 5 \quad \textmd{ and } \quad \alpha \leq \sqrt{\frac{n-2}{d-2}}
\end{equation} 
then $u$ is radially symmetric and it is given by \eqref{u_radial} (up to a translation in the case $a=b=0$). 
\end{theorem}

We mention that condition \eqref{alpha_reduced} on $\alpha$ already appeared in \cite{CC} and a geometric interpretation follows by a conformal reformulation of the problem, as it was noticed in \cite{DGZ}. Indeed, it may be convenient to write \eqref{CKN} and \eqref{eq_CKN} in a Riemannian setting, where we define the conformal metric
\begin{equation} \label{g_def}
g_{ij} = |x|^{2(\alpha-1)} \delta_{ij} \,,
\end{equation}
which permits to write \eqref{CKN} as a weighted Sobolev inequality on the manifold $(\R^d, g)$: 
\begin{equation} \label{CKN_g}
\bigg(\int_{\R^d}|u|^pd\mu(x)\bigg)^\frac{2}{p}\le C_{a,b}\int_{\R^d}|\nabla_{{g}} u|_g^2d\mu(x) \,,
\end{equation}
with $d\mu(x)=|x|^{-bp}dx$. It can be shown that $(\R^d, g, \mu)$ satisfies a $CD(0,n)$ condition if and only if $\alpha$ satisfies \eqref{alpha_reduced}. We mention that the condition \eqref{alpha_optimal} corresponds to a weaker form of the $CD(0,n)$ condition and we refer to \cite{DGZ} for a more detailed discussion on this topic.

In order to prove Theorem \ref{thm_main} we take inspiration from \cite{CFP,CMR,OU, Vetois} and we employ a $P$-function approach, as we are going to explain in the following. After the change of metric \eqref{g_def}, \eqref{eq_CKN} can be written as
\begin{equation*}
\Cl u=u^\frac{n+2}{n-2} \,,
\end{equation*}
where
\begin{equation*}
\Cl u : = -e^f \textrm{div}_{ g}\bigg(e^{-f}\nabla_{{g}}u\bigg) \,, \quad \textmd{with } e^{-f}:=|x|^{\alpha(n-d)} \,.
\end{equation*}
By letting $v=u^{-\frac{2}{n-2}}$ we obtain that 
\begin{equation}\label{eq_v}
\Cl v=\frac{2}{n-2}\frac{1}{v}+\frac{n}{2}\frac{|\nabla_{{g}} v|^2_{ g}}{v}.
\end{equation}
The advantage of considering $v$ instead of $u$ relies upon the fact that Theorem \ref{thm_main} is proved once we show that the Hessian of $v$ (denoted by $H_g v$) is a multiple of $g$, and this leads to an explicit representation of $v$ which corresponds to the fact that $u$ is given by \eqref{u_radial}. Moreover, in this case, the RHS of \eqref{eq_v} is constant and may serve as a $P$-function. Indeed, by letting 
\begin{equation} \label{P_def}
P=\frac{2}{n-2}\frac{1}{v}+\frac{n}{2}\frac{|\nabla_{{g}} v|^2_{ g}}{v}
\end{equation}
we prove that $P$ is a subsolution of a suitable elliptic PDE, and this permits us to prove the rigidity result by showing that $P$ is constant as in some classical Liouville Theorem (see for instance \cite[Proposition 8.2]{Farina}).

Theorem \ref{thm_main} is a particular case of a more general theorem where we consider \eqref{eq_CKN} defined in a convex cone. The study of Sobolev inequalities in convex cones was initiated by Lions, Pacella and Tricarico in \cite{LPT} and then more general results were obtained by the first author together with Figalli and Roncoroni in \cite{CFR} (under the finite energy assumption on the solution). In this paper, we obtain a classification of critical points under the assumption \eqref{alpha_reduced} and by assuming that the cone is convex. 

More precisely, let $A\subset \mathbb{S}^{d-1}$ be a smooth and open domain in the $(d-1)$-dimensional sphere and consider the cone $\Sigma$ defined by
$$
\Sigma = \{x \in \R^d \,:\ x=t \xi \ \textmd{for some } t \in (0,+\infty) \textmd{ and } \xi \in A \} \,.
$$ 
We notice that it may be convenient to write the cone $\Sigma$ as
$$
\Sigma=\mathbb{R}^k\times\mathcal{C} \,,
$$
where $k\in\{0,\dots,d\}$ and $\mathcal C\subset\mathbb{R}^{d-k}$ is a convex cone centred at $\bold O \in \mathbb{R}^{d-k}$ (the origin of $\mathbb{R}^{d-k}$) which does not contain a line, and we say that $x_0 \in \mathbb R^k \times \mathcal C $ is a \emph{vertex} of $\Sigma$. In the case $k=d$ we have $\Sigma = \mathbb R^d$ and a vertex of $\Sigma$ is any point in $\mathbb R^d$. 

Hence, in the case of a cone $\Sigma$, we consider positive functions in 
\begin{equation*}
\mathcal{D}^{a,b}_{loc} (\overline\Sigma) : = \left\{ u \in L^p_{loc}(\overline \Sigma,|x|^{-b}dx) \textmd{ such that  } |x|^{-a} |Du| \in L^2_{loc}(\overline\Sigma,dx) \right\} \,,
\end{equation*} 
and we study the following problem
\begin{equation} \label{eq_CKN_cone}
\begin{cases}
\divergence \left( |x|^{-2a} Du \right) + |x|^{-bp} u^{p-1} = 0 & \textmd{ in } \Sigma \\
u_\nu = 0  & \textmd{ on } \partial \Sigma \,.
\end{cases}
\end{equation}
We mention that the homogeneous Neumann condition on $\partial \Sigma$ is the natural one which appears when one writes the Euler-Lagrange equation associated to CKN inequalities in cones. 

Our second main result is the following.

\begin{theorem} \label{thm_main2}
Let $\Sigma \subset \R^d$ be a convex cone, and let $u \in \mathcal{D}^{a,b}_{loc} (\Sigma)$ be a positive solution to \eqref{eq_CKN_cone}. Let $\alpha$ and $n$ be given by \eqref{alpha_def} and \eqref{n_def}, respectively. If $n$ and $\alpha$ satisfy \eqref{alpha_reduced} then $u$ is radially symmetric and is given by \begin{equation} \label{u_radial_cone}
u(x) = \left( 1+ |x-x_0|^{(p-2)(a_c-a)}\right)^{-\frac{2}{p-2}}  \,,
\end{equation}
where $x_0$ is a vertex of $\Sigma$, up to rescalings.
\end{theorem}

The proof of Theorem \ref{thm_main2} follows the same argument as the one of Theorem \ref{thm_main}. The presence of the boundary condition can be tackled in some standard way and does not add substantial difficulties to the proof. However, it is worth noticing that the convexity of the cone appears in two different ways during the proof.

The first one is when we prove that $P$ is a subsolution of an elliptic PDE, since an extra term appears in the calculations. This quantity can be written in terms of the second fundamental form of $\partial \Sigma$, which is pointwise nonnegative thanks to the convexity of the cone, and then it can be easily tackled. 

The convexity of the cone also appears in another step of the proof, in particular when we use a Poincar\'e inequality on the set $A \subset \mathbb{S}^{d-1}$. In this case, it is crucial that the first nontrivial Neumann eigenvalue $\lambda_1(A)$ of the Laplace-Beltrami operator on $A$ is greater or equal than $d-1$ (see \cite{DEL}). This is true both for $\mathbb{S}^{d-1}$ (which corresponds to Theorem \ref{thm_main}) and whenever $A$ is a convex domain on the sphere (see \cite[Theorem 4.3]{Escobar}).

We mention that, in the case of Sobolev inequality (i.e. $a=b=0$),  if $\lambda_1(A) < d-1$ then it was proved in \cite{CPP} that symmetry breaking occurs. It is an open question whether or not rigidity results can be obtained in non-convex cones which have $\lambda_1(A) > d-1$, also under the assumption of finite energy. This result would be new also for the classical Sobolev inequality.

\medskip

The paper is organized as follows. In Section \ref{section2} we introduce the conformal metric and prove some preliminary results. Section \ref{section3} is devoted to the $P$-function and its properties. In Section \ref{section4} we prove Theorem \ref{thm_main}. The proof of Theorem \ref{thm_main2} is sketched in Section \ref{section5}.

\section{The conformal formulation} \label{section2}
Let $g$ be the metric given by \eqref{g_def}, i.e.
\begin{equation*}
g_{ij} = |x|^{2(\alpha-1)} \delta_{ij} \,,
\end{equation*}
where $\alpha$ is given by \eqref{alpha_def} and $\delta$ is the Euclidean metric. Then \eqref{CKN} can be written as 
\begin{equation}\label{CKN_g2}
\left(\int_{\R^d}|u|^p e^{-f} dV_g\right)^\frac{2}{p}\le C_{a,b}\int_{\R^d}|\nabla_{{g}} u|_g^2e^{-f} dV_g \,,
\end{equation}
where 
\begin{equation} \label{dVg}
dV_{ g}=|x|^{d(\alpha-1)}dx \,,
\end{equation}
and  
\begin{equation} \label{e-f_def}
e^{-f} =|x|^{\alpha(n-d)} \,,
\end{equation}
with $n$ given by \eqref{n_def}. In order to emphasize that we are considering the Riemannian manifold $(R^d,g)$, it will be convenient to rename the set $\mathcal D^{a,b}(\R^d)$ as $\mathcal D^g(\R^d)$, where
\begin{equation} \label{Dg} 
\mathcal D^g(\R^d) = \left\{ u: \R^d \to \R :\  \int_{\R^d}|u|^p e^{-f} dV_g < + \infty \ \textmd{ and }\  \int_{\R^d}|\nabla_{{g}} u|_g^2e^{-f} dV_g < +\infty \right\} \,,
\end{equation}
and we set $\mathcal D_{loc}^g(\R^d)$ accordingly. Correspondingly, \eqref{eq_CKN} becomes
\begin{equation} \label{eq_g_1}
\textrm{div}_{g}\bigg(e^{-f} \nabla_{{g}}u\bigg)+e^{-f} u^\frac{n+2}{n-2}=0
\end{equation}
in $\R^d$, with $u \in \mathcal D_{loc}^g(\R^d)$. It will be convenient to set 
\begin{equation} \label{L_def}
\Cl u = e^f \textrm{div}_{g}\bigg(e^{-f} \nabla_{{g}}u\bigg) \,,
\end{equation}
and hence \eqref{eq_g_1} can be written as
\begin{equation} \label{eq_g_L}
\Cl u + u^{\frac{n+2}{n-2}}=0 \quad  
\end{equation}
in $\R^d$. We notice that one can also write $\Cl u$ as follows
\begin{equation} \label{Lv_expand} 
\Cl u = \Delta_g u-\nabla_g f\cdot_g\nabla_g u \,.
\end{equation}
Hence, by a weak solution of \eqref{eq_g_L} we mean a function $u \in \mathcal D_{loc}^g(\R^d)$ such that
$$
- \int_{\R^d} \nabla_{g} u \cdot_{g}\nabla_{g} \varphi \, e^{-f} dV_g + \int_{\R^d} u^{\frac{n+2}{n-2}} \varphi  \, e^{-f} dV_g = 0 \,,
$$
for any $\varphi \in C_c^\infty (\R^d)$, where $\cdot_g$ denotes the inner product in $(\R^d,g)$. By standard elliptic regularity, solutions of \eqref{eq_CKN} are smooth in $\R^d \setminus \{O\}$ and it can be shown that they are locally bounded (see for instance \cite{CC}). The same regularity results apply to solutions to \eqref{eq_g_L}.

We conclude this section by giving some information on the weighted manifold $(\R^d, g, e^{-f})$. We first write the Ricci curvature $Ric_g$ and the Hessian $H_g$ of a function by using that $g$ is conformal to the Euclidean metric $\delta$, in particular by writing
$$
{g}=|x|^{2(\alpha-1)}\delta=e^{2\phi}\delta \,,\quad \textmd{with} \quad  \phi:=(\alpha-1)\log |x|. 
$$
In this setting we recall that for conformal metrics $g$ and $\delta$ we have (see for instance \cite[Formulas (2.68) and (2.74)]{CatinoMastrolia}) 
\begin{equation} \label{Ric_conformal}
Ric_{{g}}=Ric_\delta-(d-2)[\phi_{ij}-\phi_i\phi_j]_{i,j=1}^d-(\Delta \phi+(d-2)|\nabla \phi|^2)\delta \,,
\end{equation}
with $Ric_\delta=0$ since it is the Euclidean metric, and for any $C^2$ function $h$ we have
%\begin{equation} \label{Hessian_conformal}
%H_{{g}}h=H_\delta h+(\nabla\phi \cdot \nabla h )\delta-2[\phi_i h_j]_{i,j=1}^d
%\end{equation}
%\todo{ATTENZIONE -- FORMULA GIUSTA E':}
\begin{equation} \label{Hessian_conformal}
H_{{g}}h=H_\delta h+(\nabla\phi \cdot \nabla h )\delta- (d \phi \otimes d h + d h \otimes d \phi) \,,
\end{equation}
%\todo{CHECK SE CAMBIA QUALCOSA -- perche' l'hai scritta con il 2?}
where $H_\delta h$ is the standard Hessian of $h$ in $\R^d$. Hence, straightforward calculations show that 
\begin{equation} \label{ricci}
(Ric_g)_{ij} = \frac{(d-2)(1-\alpha^2)}{|x|^2}   \left(\delta_{ij} - \frac{x_i x_j}{|x|^2} \right) 
\end{equation}
and, by recalling that $f$ in \eqref{e-f_def} is given by
\begin{equation} \label{f_expl}
f=-\alpha(n-d)\log |x| \,,
\end{equation}
we have
\begin{equation} \label{Hess_f}
(H_{{g}}f)_{ij}=-\alpha^2\frac{n-d}{|x|^2}\bigg(\delta_{ij}-2\frac{x_i x_j}{|x|^2}\bigg) \,.
\end{equation}
We notice that, since $d \geq 2$, then we have that $Ric_g \geq 0$ for $\alpha \leq 1$. Moreover, we have the following lemma.

\begin{lemma} \label{lemma_k}
Let $v \in C^0(\R^d) \cap C^3(\R^d \setminus \{O\})$ and set  
\begin{equation} \label{k_def}
k[v]:=\big{|}H_{{g}} v\big{|}_g^2 -\frac{(\Cl v)^2}{n} +Ric_{{g}}(\nabla_{{g}} v,\nabla_{{g}} v)+H_{{g}}f(\nabla_{{g}} v,\nabla_{{g}} v) \,.
\end{equation}
Then 
%\begin{equation} \label{kgeqW}
%k[v]\geq |H_{{g}} v-\frac{\Delta_g v}{d}g|_{{g}}^2+W_f[v] \,,
%\end{equation}
\begin{equation} \label{k_nuova}
k[v] = \Big{|}H_{{g}} v-\frac{\Delta_g v}{d}g \Big{|}_{{g}}^2+\frac{n-d}{nd} \left(\Delta_g v -\alpha d  \frac{\nabla_g v \cdot x}{|x|^2} \right)^2 + W_f[v] \,, 
\end{equation}
where
\begin{equation} \label{Wf_def}
W_f[v]= \frac{d-2-\alpha^2(n-2)}{|x|^2}\bigg(|\nabla_{ g}v|^2-\frac{(\nabla_{ g}v\cdot x)^2}{|x|^2}\bigg) \,.
\end{equation}
Moreover, if 
\begin{equation}\label{condalpha}
\alpha^2\le \frac{d-2}{n-2}
\end{equation}
then 
\begin{equation} \label{kgeq0}
k[v]\ge 0
\end{equation}
and 
\begin{equation} \label{v_radial}
k[v] = 0 \quad \textmd{ if and only if } \quad v(x)=
\begin{cases}
c_1|x-x_0|^{2\alpha} + c_2 & \textmd{ if } a=b=0 \,, \\
c_1|x|^{2\alpha} + c_2 & \textmd{ otherwise} \,, 
\end{cases}
\end{equation}
for any $ x\in \R^d$ for some $c_1,c_2\in\R$ and $x_0 \in \R^d$.
\end{lemma}
\begin{proof}
We first notice that, from the definition of $\Cl$ \eqref{L_def}, it follows:
\begin{equation} \label{eq_treno}
\begin{split}
|H_{{g}} v|_{{g}}^2-\frac{(\Cl v)^2}{n}& =|H_{{g}} v|_{{g}}^2-\frac{(\Delta_g v)^2}{d}+\frac{(\Delta_g v)^2}{d}-\frac{1}{n}\bigg((\Delta_g v)^2-2\Delta_g v(\nabla_g v \cdot_g\nabla_g f)+(\nabla_g v \cdot_g\nabla_g f)^2\bigg) \\
& =\Big{|}H_{{g}} v-\frac{\Delta_g v}{d}g \Big{|}_{{g}}^2+\frac{n-d}{nd} (\Delta_g v)^2 +\frac{2}{n}\Delta_g v(\nabla_g v \cdot_g\nabla_g f)-\frac{1}{n}(\nabla_g v \cdot_g\nabla_g f)^2.
\end{split}
\end{equation}
In the following, we assume that $n\neq d$, since otherwise $f=0$ and \eqref{k_nuova} immediately follows. Hence from \eqref{k_def}, \eqref{ricci} and \eqref{Hess_f}, we obtain that $k[v]$ and can be written as 
\begin{multline*}
k[v]= \Big{|}H_{{g}} v-\frac{\Delta_g v}{d}g \Big{|}_{{g}}^2+\frac{n-d}{nd} (\Delta_g v)^2 +\frac{2}{n}\Delta_g v(\nabla_g v \cdot_g\nabla_g f) + \frac{d}{n(n-d)}(\nabla_g v \cdot_g\nabla_g f)^2 \\ 
- \frac{1}{n-d} (\nabla_g v \cdot_g\nabla_g f)^2 + Ric_{{g}}(\nabla_{{g}} v,\nabla_{{g}} v)+H_{{g}}f(\nabla_{{g}} v,\nabla_{{g}} v)
\end{multline*}
and from \eqref{f_expl} we find
\begin{multline*}
k[v] = \Big{|}H_{{g}} v-\frac{\Delta_g v}{d}g \Big{|}_{{g}}^2+\frac{n-d}{nd} \left(\Delta_g v + \frac{d}{n-d} \nabla_g v \cdot_g\nabla_g f \right)^2 \\
- \frac{\alpha^2(n-d)}{|x|^2} \left(\nabla_g v \cdot \frac{x}{|x|}\right)^2 + Ric_{{g}}(\nabla_{{g}} v,\nabla_{{g}} v)+H_{{g}}f(\nabla_{{g}} v,\nabla_{{g}} v) \,. 
\end{multline*}
From \eqref{ricci} and \eqref{Hess_f} we find \eqref{k_nuova}; from \eqref{condalpha} and Cauchy-Schwarz inequality we obtain that $W_f [v]\geq 0$ and then \eqref{kgeq0}.
%
%
%\vspace{3em}
%
%Young's inequality yields
%$$
%\frac{2}{n}\Delta_g v(\nabla_g v \cdot_g\nabla_g f)\ge -\frac{n-d}{nd}(\Delta_g v)^2-\frac{d}{n(n-d)}(\nabla_g v \cdot_g\nabla_g f)^2
%$$
%and then from \eqref{eq_treno} we obtain
%$$
%|H_{{g}} v|_{{g}}^2-\frac{(\Cl v)^2}{n}\ge |H_{{g}} v-\frac{\Delta_g v}{d}g|_{{g}}^2-\frac{1}{n-d} (\nabla_g v \cdot_g\nabla_g f)^2.
%$$
%Hence from \eqref{k_def} we find  
%\begin{equation} \label{wait}
%k[v]\ge  \Big{|}H_{{g}} v-\frac{\Delta_g v}{d}g \Big{|}_{{g}}^2-\frac{1}{n-d} (\nabla_g v \cdot_g\nabla_g f)^2+Ric_{{g}}(\nabla_{{g}} v,\nabla_{{g}} v)+H_{{g}}f(\nabla_{{g}} v,\nabla_{{g}} v) \,,
%\end{equation}
%and, by using \eqref{ricci} and \eqref{Hess_f}, we obtain \eqref{kgeqW}. Moreover, \eqref{kgeq0} immediately follows from \eqref{kgeqW}, \eqref{Wf_def}, \eqref{condalpha} and by using Cauchy-Schwarz inequality.

Now we assume that $k[v]=0$. From \eqref{k_nuova} we have
\beq\label{comp41}
|H_{{g}} v-\frac{\Delta_g v}{d}g|_{{g}}^2=0 \,,
\eeq
%
%-----------
%
%\eqref{kgeq0} 
%$$
%\Delta_g v -\alpha d  \frac{\nabla_g v \cdot x}{|x|^2} = 0\,,
%$$
%and 
%\beq\label{comp42}
%\alpha^2 < \frac{d-2}{n-2} \quad \textmd{ and } \quad |\nabla_{ g}v|^2= \frac{(\nabla_{ g}v\cdot x)^2}{|x|^2} \,.
%\eeq
%or \todo{c'e' anche questo caso}
%$$
%\alpha^2 = \frac{d-2}{n-2} \,.
%$$
%\todo{pero' non mi sembra che si usi la \eqref{comp42}}
%
%--------------
%
By using \eqref{Hessian_conformal}, from \eqref{comp41} it follows 
\beq\label{comp44}
0\equiv H_{{g}} v-\frac{\Delta_g v}{d}g=H v-\frac{\Delta v}{d}\delta-2[\phi_i v_j]_{i,j=1}^d+\frac{2}{d}(\nabla \phi\cdot\nabla v)\delta \,,
\eeq
where we also used that
$$
e^{2\phi}\Delta_g v = \Delta v + (d-2) \nabla \phi\cdot\nabla v  \,.
$$
By multiplying by $e^{-2\phi}$, we notice that \eqref{comp44} can be written as
\begin{equation} \label{bicocca}
(e^{-2\phi}v_i)_j=\frac{1}{d}(e^{-2\phi}v_k)_k\delta_{ij}
\end{equation}
for $i,j \in \{1,\ldots,d\}$. If $i \neq j$ then \eqref{bicocca} 
\begin{equation} \label{esatto}
(e^{-2\phi}v_i)_j= 0  \quad \textmd{ if } i \neq j \,.
\end{equation}
If $i=j$, since $v$ and $\phi$ are smooth outside the origin, we can differentiate and obtain that 
$$
(e^{-2\phi}v_i)_{i\ell} = \frac{1}{d}(e^{-2\phi}v_k)_{k\ell}
$$
for $i,\ell =1,\ldots,d$. From \eqref{esatto} we find
$$
(e^{-2\phi}v_i)_{ii} = 0 \quad \textmd{for } i = 1,\ldots,d \,,
$$
which implies
\begin{equation} \label{pranzo}
(e^{-2\phi}v_i)_{i} = C_1^i \quad \textmd{for } i = 1,\ldots,d \,.
\end{equation}
for some constant vector $C_1=(C_1^1,\ldots,C_1^d)$. Hence from \eqref{pranzo} and \eqref{esatto} we have
\begin{equation} \label{grad_v}
v_i = |x|^{2(\alpha-1)}(C_1^i x_i + C_2),
\end{equation}
for some constant vector $C_2$. In order to complete the proof, we need to distinguish between three cases, according to the value of $\alpha$.

Since $k[v]=0$ and if we assume 
$$
\alpha^2 < \frac{d-2}{n-2}
$$
then we have $W_f[v]=0$, which yields 
\begin{equation*}
|\nabla_{ g}v|^2= \frac{(\nabla_{ g}v\cdot x)^2}{|x|^2} 
\end{equation*}
and \eqref{grad_v} implies that $C_2=0$.

If 
$$
\alpha^2 = \frac{d-2}{n-2} \quad \textmd{ and } n \neq d
$$
then $k[v]=0$ implies that 
$$
\Delta_g v = \alpha d  \frac{\nabla_g v \cdot x}{|x|^2}  \,.
$$
Since 
$$
\Delta_g v = d e^{-2\phi}  \phi_i v_i +  \sum_{i=1}^d (e^{-2\phi}v_i)_{i}
$$
from \eqref{pranzo} and \eqref{grad_v} we obtain again that $C_2=0$.

Finally, we assume that $n=d$ and $\alpha=1$, i.e. $a=b=0$ as it follows from \eqref{alpha_def} and \eqref{n_def}. From \eqref{grad_v} we immediately obtain that $v(x)=c_1 |x-x_0|^2 + c_2$, for some constants $c_1,c_2$ and for some $x_0 \in \R^d$.
\end{proof}

\section{The p-function} \label{section3}
In this section, we introduce an auxiliary function, the so-called $P$-function and we prove some preliminary results.
Let $u$ be a solution to \eqref{eq_g_L} and set
\begin{equation} \label{v_def}
v=u^{-\frac{2}{n-2}} \,,
\end{equation}
it is readily seen that $v$ satisfies
\begin{equation}\label{eq2}
\Cl v=\frac{2}{n-2}\frac{1}{v}+\frac{n}{2}\frac{|\nabla_{{g}} v|^2_{ g}}{v} \,,
\end{equation}
where $\Cl$ is given by \eqref{L_def}.
Since $u$ is positive, smooth in $\R^d \setminus \{O\}$ and locally bounded in $\R^d$, then the same holds for $v$. Moreover, $u\in \mathcal D_{loc}^g(\R^d)$ implies that $v^{-\frac{n-2}{2}} \in \mathcal D_{loc}^g(\R^d)$, i.e.
\begin{equation} \label{v_Dg}
\int_{B_\rho} v^{-\frac{n-2}{2}p} e^{-f} dV_g < +\infty \quad \textmd{ and } \quad  \int_{B_\rho} v^{-n} |\nabla_{{g}} v|_g^2  e^{-f} dV_g < +\infty 
\end{equation}
for any ball of radius $\rho>0$.

The right-hand side of \eqref{eq2} is what we call $P$-function, namely we set
\begin{equation}\label{P-fun}
P=\frac{2}{n-2}\frac{1}{v}+\frac{n}{2}\frac{|\nabla_{{g}} v|^2_{ g}}{v} \,.
\end{equation}

Thus $P=\Cl v$ and our main goal in the proof of Theorem \ref{thm_main} will be to show that $P$ is constant. This is essentially the approach proposed in \cite{CFP} for the critical Laplace equation in completes manifolds with nonnegative Ricci curvature. In our setting, the manifold is weighted and it is not complete, for this reason, some further arguments are needed. 

The crucial property of $P$ is that $P$ is a subsolution to a suitable elliptic equation, as it is proved in the following lemma.

\begin{lemma} \label{lemma_P}
Let $v$ be a solution to \eqref{eq2} satisfying \eqref{v_Dg} and let $P$ be defined by \eqref{P-fun}.  Then $P$ satisfies
\begin{equation}\label{eqfond}
e^f \divergence_{g}\Big(e^{-f} v^{2-n}\nabla_{ g}P\Big)=nv^{1-n}k[v] \quad \textmd{ in } \R^d \setminus \{O\} \,,
\end{equation}
where $\Cl$ is given by \eqref{L_def}.
\end{lemma}

\begin{proof}
We recall that $v$ is smooth in $\R^d \setminus \{O\}$ and hence we can verify \eqref{eqfond} pointwise in $\R^d \setminus \{O\}$. From Lemma \ref{lemma_app} we have
\begin{equation} \label{2.18}
e^f \divergence \bigg( e^{-f} v^{1-n}\bigg( \frac{n}{2}\nabla_{{g}}|\nabla_{{g}} v|^2_{{g}} - n \Cl v\nabla_{{g}} v + (n-1)P \nabla_{{g}} v\bigg)\bigg) \\
=n v^{1-n}k[v] + v (\Cl v- P ) \Cl(v^{1-n}),
\end{equation}
where $k[v]$ is given by \eqref{k_def}. From \eqref{eq2} and \eqref{2.18} we get
\begin{equation}\label{2}
e^{f}\textrm{div}_{{g}}\bigg(e^{-f}v^{1-n}\bigg(\frac{1}{2}\nabla_{{g}}|\nabla_{{g}} v|^2_{{g}} - \frac{1}{n}\Cl v\nabla_{{g}} v\bigg)\bigg)=v^{1-n}k[v].
\end{equation}
Since $\Cl v =P$ and from \eqref{P-fun}, we have that
$$
\nabla_g P =   \frac{n}{2v}\nabla_{{g}}|\nabla_{{g}} v|^2_{{g}}  -\frac{1}{v} P \nabla_g v = \frac{n}{v} \left( \frac{1}{2}\nabla_{{g}}|\nabla_{{g}} v|^2_{{g}} -\frac{1}{n} \Cl v \nabla_g v \right)
$$
which, together with \eqref{2}, implies \eqref{eqfond}.
\end{proof}

Lemma \ref{lemma_P} is preliminary to the following lemma, which provides the crucial differential inequality that we are going to use.

\begin{lemma} \label{lemma_Pt}
Let $v$ be a solution to \eqref{eq2} satisfying \eqref{v_Dg} and let $P$ be defined by \eqref{P-fun}. Then 
\begin{equation} \label{eqPt}
e^f \divergence_{g}\Big(e^{-f} P^{t-1}  v^{2-n}\nabla_{ g}P\Big)  \geq \left(t-\frac12\right) P^{t-2} v^{2-n}|\nabla_{ g}P|_g^2 + nP^{t-1}v^{1-n} W_f[v]    
\end{equation}
in $\R^d \setminus \{O\}$ and for any $t \in \R$, where $\Cl$ is given by \eqref{L_def}.
\end{lemma}

\begin{proof}

We first notice that $v\in C^\infty(\R^d\setminus\{O\})$ and then the following computations hold pointwise in $\R^d\setminus\{O\}$. 
From the definition of the $P$-function \eqref{P-fun} we have that 
$$
\nabla_{{g}}P=\frac{n}{v}\bigg(H_{{g}} v-\frac{\Cl v}{n} g \bigg) \nabla_{{g}} v.
$$
By applying Cauchy-Schwarz inequality we get
\begin{equation*}
|\nabla_{{g}} P|_{{g}}^2\le\frac{n^2}{v^2}|H_{{g}}v-\frac{\Cl v}{n}{g}|_{{g}}^2|\nabla_{{g}} v|_{{g}}^2
\end{equation*}
and, since $P>0$, we have
\begin{equation}\label{comp1}
P^{t-2}v^{2-n}|\nabla_{{g}} P|_{{g}}^2\le P^{t-2}v^{2-n}\frac{n^2}{v^2}|H_{{g}}v-\frac{\Cl v}{n}{g}|_{{g}}^2|\nabla_{{g}} v|_{{g}}^2= P^{t-1}v^{1-n}n^2|H_{{g}}v-\frac{\Cl v}{n} g|_{{g}}^2\frac{|\nabla_{{g}} v|_{{g}}^2}{vP}.
\end{equation}
From the definition of $P$ \eqref{P-fun} we have
\begin{equation}\label{comp2}
\frac{|\nabla_{{g}} v|_{{g}}^2}{vP}\le\frac{2}{n},
\end{equation}
and then \eqref{comp1} and \eqref{comp2} yield
\beq\label{comp3}
P^{t-2}v^{2-n}|\nabla_{{g}} P|_{{g}}^2\le P^{t-1}v^{1-n}2n|H_{{g}}v-\frac{\Cl v}{n} g|_{{g}}^2.
\eeq
From \eqref{Lv_expand} we get
\begin{equation*}
\begin{split}
|H_{{g}}v-\frac{\Cl v}{n} g|_{{g}}^2 & =|H_{{g}}v|_{{g}}^2+\frac{d}{n^2}(\Delta_g v-\nabla_g f\cdot_g\nabla_g v)^2-\frac{2}{n}(\Delta_g v-\nabla_g f\cdot_g\nabla_g v) \Delta_gv \\
 & =|H_{{g}}v|_{{g}}^2+\bigg(\frac{d}{n^2}-\frac{2}{n}\bigg)(\Delta_g v)^2+\frac{d}{n^2}(\nabla_g f\cdot_g\nabla_g v)^2+\bigg(\frac{2}{n}-\frac{2d}{n^2}\bigg)(\nabla_g f\cdot_g\nabla_g v) \Delta_gv \\ 
 & =|H_{{g}}v- \frac{\Delta_g v}{d } g|_{{g}}^2+\frac{(n-d)^2}{dn^2}(\Delta_g v)^2+\frac{d}{n^2}(\nabla_g f\cdot_g\nabla_g v)^2+ 2 \frac{n-d}{n^2} (\nabla_g f\cdot_g\nabla_g v) \Delta_gv \,,
\end{split}
\end{equation*}
and hence
\begin{equation*}
|H_{{g}}v-\frac{\Cl v}{n} g|_{{g}}^2 = |H_{{g}}v- \frac{\Delta_g v}{d } g|_{{g}}^2+ \frac{1}{dn^2} \left[ (n-d)\Delta_g v +d (\nabla_g f\cdot_g\nabla_g v) \right]^2 \,.
\end{equation*}	
From \eqref{comp3} we have
\begin{equation*} 
P^{t-2}v^{2-n}|\nabla_{{g}} P|_{{g}}^2\le P^{t-1}v^{1-n}2n \left\{ |H_{{g}}v- \frac{\Delta_g v}{d } g|_{{g}}^2+ \frac{1}{dn^2} \left[ (n-d)\Delta_g v +d (\nabla_g f\cdot_g\nabla_g v)  \right]^2 \right\} \,,
\end{equation*}
and, by using \eqref{k_nuova}, we obtain
\begin{equation} \label{starb}
P^{t-2}v^{2-n}|\nabla_{{g}} P|_{{g}}^2 + 2nP^{t-1}v^{1-n} W_f[v] \le 2n P^{t-1}v^{1-n} k[v] \,.
\end{equation}

Since 
\begin{equation*}
e^f \divergence_{g}\Big(e^{-f}P^{t-1}  v^{2-n}\nabla_{ g}P\Big)  =nP^{t-1} v^{1-n}k[v] + (t-1) P^{t-2} v^{2-n}|\nabla_{ g}P|_g^2 \,,
\end{equation*}
from \eqref{eqfond} and \eqref{starb} we obtain \eqref{eqPt}.
\end{proof}

We shall use Lemmas \ref{lemma_P} and \ref{lemma_Pt} by testing \eqref{eqfond} with some cut off function which has compact support in $B_R^g \setminus \overline B_r^g  \subset \R^d \setminus \{O\}$ (see Proposition \ref{mainineq}). The lack of regularity at the origin implies that we need to have precise asymptotic information on some relevant quantities as $r\to0$. This is the content of the following proposition, which was proved in \cite[Section 8]{DEL} and we write here below by using the notation adopted in this paper. 
%
%We consider the spherical coordinates:
%$$
%\rho=|x|\quad\text{and}\quad \omega=\frac{x}{|x|}
%$$
%and we define: \todo{perche' introduci $s$?? A cosa ci serve??}
%\beq\label{defs}
%s(\rho,\omega):=(n-1)u^{-\frac{2}{n-2}}(\rho,\omega)=(n-1)v(\rho,\omega)
%\eeq
%Then we have
%$$
%|\nabla v|^2=\frac{1}{(n-1)^2}|\nabla s|^2=\frac{1}{(n-1)^2}|s'|^2+\frac{1}{(n-1)^2}\frac{1}{r^2}|\nabla_\omega s|^2
%$$
%where $s'=\frac{\partial s}{\partial \rho}$ and $\nabla_\omega$ denotes the gradient with respect to the angular variable $\omega\in\S^{d-1}$. On $\R^d$ we consider the measure 
%$$
%d\mu:=\rho^{n-1}d\rho d\omega.
%$$
%
%\todo{Questo va riscritto in $v$, cioe' come lo usiamo}
%
%\begin{proposition}{(\cite{DEL},Proposition 8.2, (i) and (ii))}\label{DEL}
%
%Let $v\in $\todo{Spazio} be a weak solution of  \eqref{eq2}, and $s$ defined by \eqref{defs}. 
%
%Then, if $\alpha\le\sqrt{\frac{d-1}{n-1}} $ and $\rho\to 0$
%$$
%\int_{\S^{d-1}}|s'|^2d\omega=O(1)
%$$
%and
%$$
%\int_{\S^{d-1}} |\nabla_\omega s|^2 d\omega=O(\rho^2).
%$$
%\end{proposition}
%

\begin{proposition} \label{DEL}
Let $v$ be a solution to \eqref{eq2} satisfying \eqref{v_Dg} and assume that 
$$
\alpha\le\sqrt{\frac{d-1}{n-1}}\,. 
$$
Then
$$
\int_{B_{r}^g} e^{-f}|\nabla_{ g} v|^2_{ g} dV_g \le Cr^{n-2+\frac{2}{\alpha}} \,.
$$
In particular, for $n>2$ and $\alpha \leq 1$ we have
$$
\int_{B_{r}^g} e^{-f}|\nabla_{ g} v|^2_{ g} dV_g = o(r^2)
$$
as $r\to 0$.
\end{proposition}

\begin{proof}
This proposition directly follows from \cite[Proposition 8.2, (i) and (ii)]{DEL}. Indeed one can check that, up to a multiplicative constant, the function $\bold p$ in \cite{DEL} is equal to $v$, more precisely
\begin{equation} \label{8.53}
\bold p = (n-1) v \,.
\end{equation}
Moreover, in terms of the spherical coordinates
$$
\rho=|x|\quad\text{and}\quad \omega=\frac{x}{|x|} \,,
$$  
\cite[Proposition 8.2, (i) and (ii)]{DEL} gives the existence of a constant $C$ (not depending on $\rho$) such that
\begin{equation} \label{8.54}
\int_{\mathbb{S}^{d-1}} |\partial_\rho \bold p(\rho, \omega)|^2 d\sigma \leq C
\quad \textmd{ and } \quad
\int_{\mathbb{S}^{d-1}} |\nabla_\omega \bold p(\rho, \omega)|^2 d\sigma \leq C \rho^2 \,,
\end{equation}
for as $\rho \to 0$, where $\partial_\rho$ denotes the standard Euclidean radial derivative.   
Thus, \eqref{8.53} yields
$$
\int_{B_{r}^g} e^{-f}|\nabla_{ g} v|^2_{ g} \, dV_g=\int_0^{r^{1/\alpha}} d\rho\int_{\S^{d-1}}\rho^{\alpha n-d+2-2\alpha} \bigg(\frac{1}{(n-1)^2}|\partial_\rho \bold p|^2+\frac{1}{(n-1)^2}\frac{1}{\rho^2}|\nabla_\omega \bold p|^2\bigg)\rho^{d-1} d\omega
$$
and from \eqref{8.54} we conclude.
\end{proof}

\section{Proof of Theorem \ref{thm_main}} \label{section4}
In order to prove Theorem \ref{thm_main} we need some preliminary results. To simplify the notation, in the following calculations, we shall omit the volume form $dV_g$. In the following, $B_R^g$ denotes the ball centred at the origin of radius $R$ in the metric $g$. We notice that $B_R^g$ correspond to the Euclidean ball $B_{R^\frac{1}{\alpha}}$ centered at the origin and of radius $ R^\frac{1}{\alpha}$ and that 
\begin{equation} \label{VB_R}
\int_{B_R^g} e^{-f} dV_g = C_* R^n \,,
\end{equation}
where $C_*$ depends only on $d,n$ and $\alpha$.

Given $0<r<R$, we shall refer to a cutoff function $\varphi$ with compact support in $B_{2R}^g \setminus \overline B_r^g$ as a function satisfying the following conditions 
\begin{equation}\label{varphi}
\begin{cases}
\varphi=1 & \textmd{ in }\ B_R^g \setminus B_{2r}^g  \\
\varphi=0 & \textmd{ in }\  B_r^g \textmd{ and } \R^d \setminus B_{2R}^g  \\
|\nabla_g\varphi|_g \le \frac{1}{R} & \textmd{ in }\  B_{2R}^g \setminus B_{R}^g  \\
|\nabla_g\varphi|_g \le \frac{1}{r} & \textmd{ in }\   B_{2r}^g \setminus B_{r}^g   \,.
\end{cases}
\end{equation}
By a standard approximation argument, it is clear that $\varphi$ can be taken as a test function in the weak formulation of \eqref{eq2}.

\begin{proposition}\label{mainineq}
Let $v$ be a solution to \eqref{eq2} satisfying \eqref{v_Dg} and let $P$ be defined by \eqref{P-fun}. Let $\varphi$ be given by \eqref{varphi} and assume that condition \eqref{condalpha} is fulfilled. Then
\beq\label{comp26}
\frac{1}{\theta^2}(t-1/2)^2 \int_{\R^d} e^{-f}P^{t-2}v^{2-n}|\nabla_{{g}} P|_{{g}}^2\varphi^\theta\le\frac{1}{R^2}\int_{B_{2R}^g\setminus B_R^g} e^{-f}\varphi^{\theta-2}v^{2-n}P^t+\frac{1}{r^2}\int_{B_{2r}^g\setminus B_{r}^g} e^{-f}\varphi^{\theta-2}v^{2-n}P^t,
\eeq
for any $t\ge\frac{1}{2}$, $0<r<R$ and $\theta\ge 2$.
\end{proposition}
\begin{proof}
We first notice that \eqref{condalpha} implies that $W_f[v] \geq 0$, and then \eqref{eqPt} yields
\begin{equation} \label{eqPt_2}
e^f \divergence_{g}\Big(e^{-f} P^{t-1}  v^{2-n}\nabla_{ g}P\Big)  \geq \left(t-\frac12\right) P^{t-2} v^{2-n}|\nabla_{ g}P|_g^2 \,,    
\end{equation}
in $\R^d \setminus \{O\}$. Since $\varphi\in C_c^\infty (B_{2R}^g\setminus\bar B_r^g)$ and $\varphi\ge 0$, from \eqref{eqPt_2} we have 
$$
\varphi^\theta e^f\textrm{div}_{{g}}\bigg(e^{-f}v^{2-n}P^{t-1}\nabla_{{g}} P\bigg)\ge(t-\frac{1}{2})P^{t-2}v^{2-n}|\nabla_{{g}} P|_{{g}}^2\varphi^\theta  \,.
$$
By multiplying by $e^{-f}$ and integrating by parts (we omit the volume form $dV_g$), we get
\begin{multline} \label{intbyparts1}
\left(t-\frac12\right) \int_{\R^d}  e^{-f}P^{t-2}v^{2-n}|\nabla_{{g}} P|_{{g}}^2\varphi^\theta  \le - \theta\int_{\R^d} \varphi^{\theta-1} e^{-f}v^{2-n}P^{t-1}\nabla_{{g}} P \cdot_g \nabla_{{g}}\varphi \\
\le \theta\bigg(\int_{\R^d} e^{-f} 	\varphi^\theta v^{2-n}P^{t-2}|\nabla_{{g}} P|_{{g}}^2\bigg)^\frac{1}{2}\bigg(\int_{\R^d} e^{-f}\varphi^{\theta-2}v^{2-n}P^t|\nabla_{{g}} \varphi|_{{g}}^2\bigg)^\frac{1}{2}, 
\end{multline}
where in the last inequality we applied H\"older's inequality. Hence, if we assume that $t \geq \frac{1}{2}$, we obtain
\begin{equation*}
\left(t-\frac12\right)^2 \int_{\R^d}  e^{-f}P^{t-2}v^{2-n}|\nabla_{{g}} P|_{{g}}^2\varphi^\theta  \leq \int_{\R^d} e^{-f}\varphi^{\theta-2}v^{2-n}P^t|\nabla_{{g}} \varphi|_{{g}}^2  \,.
\end{equation*}
As we already mentioned, by a standard approximation argument we can choose $\varphi$ satisfying \eqref{varphi} and then we immediately obtain  \eqref{comp26}.   
\end{proof}

\begin{lemma}\label{Lemma1}
Let $v$ be a positive solution to \eqref{eq2} satisfying \eqref{v_Dg}. For any $q\in\R$ and $\psi\in C^\infty_c(\R^d\setminus \{0\})$ it holds
\begin{equation} \label{eq_Lemma1}
(\frac{n}{2}+1-q)\int_{\R^d} e^{-f}v^{-q}|\nabla_{ g} v|_{ g}^2\psi+\frac{2}{n-2}\int_{\R^d} v^{-q}\psi e^{-f}= -\int_{\R^d} v^{1-q}\nabla_{ g} v \cdot_{ g}\nabla_{ g} \psi e^{-f}.
\end{equation}
\end{lemma}
\begin{proof}
From \eqref{eq2} and \eqref{P-fun}, and by integrating by parts, we get
$$
\int_{\R^d} e^{-f} P v^{1-q}\psi=\int_{\R^d} \Cl v e^{-f} v^{1-q}\psi=\int_{\R^d}  \textrm{div}_{ g}\bigg(e^{-f}\nabla_{{g}} v\bigg)v^{1-q}\psi=
$$
$$=(q-1)\int_{\R^d} v^{-q}e^{-f}|\nabla_{ g} v|_{ g}^2\psi -\int_{\R^d} v^{1-q}\nabla_{ g} v\cdot_{ g}\nabla_{ g} \psi e^{-f}.
$$
Using again  \eqref{P-fun} we get
$$
(\frac{n}{2}+1-q)\int_{\R^d} e^{-f}v^{-q}|\nabla_{ g} v|_{ g}^2\psi+\frac{2}{n-2}\int_{\R^d} v^{-q}\psi e^{-f}= -\int_{\R^d} v^{1-q}\nabla_{ g} v \cdot_{ g}\nabla_{ g} \psi e^{-f}.
$$
\end{proof}
\begin{corollary}\label{Corollary}
Let $v$ be a positive solution to \eqref{eq2} satisfying \eqref{v_Dg}. Let $q,m\in\R$, $\psi\in C^\infty_c(\R^d\setminus\{0\})$ and let $P$ be given by \eqref{P-fun}. Then it holds
\begin{multline} \label{eq_Corollary} 
(\frac{n}{2}+1-q)\int_{\R^d} e^{-f}v^{-q}P^m|\nabla_{ g} v|_{ g}^2\psi+\frac{2}{n-2}\int_{\R^d} v^{-q}P^m\psi e^{-f} \\
= -\int_{\R^d} v^{1-q}P^m\nabla_{ g} v \cdot_{ g}\nabla_{ g} \psi e^{-f}-m\int_{\R^d} v^{1-q}P^{m-1}\nabla_{ g} v \cdot_{ g}\nabla_{ g}P \psi e^{-f} \,.
\end{multline}
\end{corollary}
\begin{proof}
We recall that, thanks to the regularity properties of $v$, $P$ is smooth and locally bounded in $\R^d \setminus \{O\}$. Hence we can choose $\psi P^m$ in place of $\psi$ in Lemma \ref{Lemma1} and \eqref{eq_Lemma1} immediately gives \eqref{eq_Corollary}.
\end{proof}

\begin{lemma}\label{Lemma2}
Let $v$ be a positive solution to \eqref{eq2} satisfying \eqref{v_Dg} and let $R>0$.

(i) If $2\le q<\frac{n}{2} + 1$ then 
\begin{equation} \label{(i)}
\int_{B_{2R}^g} e^{-f} v^{-q}|\nabla_{ g} v|_{ g}^2+\int_{B_{2R}^g}e^{-f}v^{-q}\le C R^{n-q} \,.
\end{equation}

(ii) If $ 0 \leq q\le \frac{n}{2}+1$ then
\begin{equation} \label{(ii)}
\int_{B_{2R}^g}e^{-f}v^{-q}\le C R^{n-q} \,.
\end{equation}
\end{lemma}
\begin{proof}
{\bf Step 1: proof of (i)}. 
Let $\theta >q$ and, for $0<r<R$, let $\varphi$ satisfy \eqref{varphi}. By letting $\psi = \varphi^\theta$ in Lemma \ref{Lemma1}, we have
\begin{multline} \label{sotto}
(\frac{n}{2}+1-q)\int_{\R^d} e^{-f}v^{-q}|\nabla_{ g} v|_{ g}^2\varphi^\theta+\frac{2}{n-2}\int_{\R^d} e^{-f} v^{-q}\varphi^\theta \\
=-\theta\int_{\R^d} e^{-f}\varphi^{\theta-1}v^{1-q}\nabla_{ g} v \cdot_{ g}\nabla_{ g} \varphi \leq  \theta\int_{\R^d} e^{-f}\varphi^{\theta-1}v^{1-q}|\nabla_{ g} v|_g |\nabla_{ g} \varphi|_g \,,
\end{multline}
where in the last inequality we have used Cauchy-Schwarz inequality. Let $\epsilon >0$ be a (small) constant to be chosen later. By using Young's inequality we have
\begin{equation} \label{yyoung}
\varphi^{\theta-1} v^{1-q}|\nabla_{ g} v|_{ g} |\nabla_{ g}\varphi|_g \leq \epsilon\varphi^\theta v^{-q}|\nabla_{ g} v|_{ g}^2+\frac{1}{\epsilon}  |\nabla_{ g}\varphi|_g^2 \varphi^{\theta-2}v^{2-q} \,,
\end{equation}
and a further application of Young's inequality, this time with exponents $\bigg(\frac{q}{q-2},\frac{q}{2}\bigg)$ and assuming $q>2$, yields
\begin{equation} \label{young2volte}
\varphi^{\theta-1} v^{1-q}|\nabla_{ g} v|_{ g} |\nabla_{ g}\varphi|_g \leq  \epsilon\varphi^\theta v^{-q}|\nabla_{ g} v|_{ g}^2 + \epsilon \varphi^\theta v^{-q} + \frac{1}{\epsilon^{q-1}} |\nabla_{ g}\varphi|_g^q \varphi^{\theta -q} \,.
\end{equation} 
We notice that if $q=2$ then \eqref{young2volte} immediately follows from \eqref{yyoung}.

%
%and from \eqref{varphi} we obtain
%\begin{multline} \label{comp7}
%(\frac{n}{2}+1-q)\int_{\R^d} e^{-f}v^{-q}|\nabla_{ g} v|_{ g}^2\varphi^\theta+\frac{2}{n-2}\int_{\R^d} e^{-f}v^{-q}\varphi^\theta\le \theta\epsilon\int_{\R^d} e^{-f}\varphi^\theta v^{-q}|\nabla_{ g} v|_{ g}^2 \\
%+\frac{\theta}{\epsilon R^2}\int_{B_{2R}^g\setminus B_R^g} e^{-f}\varphi^{\theta-2}v^{2-q}+\frac{\theta}{\epsilon r }\int_{B_{2 r}^g\setminus B_{r}^g} e^{-f}\varphi^{\theta-2}v^{2-q}.
%\end{multline}
%We apply Young's inequality,  with exponents $\bigg(\frac{q}{q-2},\frac{q}{2}\bigg)$, and notice that 
%\beq\label{comp8}
%\frac{1}{\epsilon \rho^2}\int_{B_{2\rho}\setminus B_\rho} e^{-f}\eta^{\theta-2}v^{2-q}\le \epsilon\int_{\R^d} e^{-f}\eta^{\theta}v^{-q}+\frac{1}{\epsilon^{q-1} \rho^q}\int_{B_{2\rho}\setminus B_\rho} e^{-f}\eta^{\theta-q}.
%\eeq
%\beq\label{comp28}
%\frac{\rho^2}{\epsilon }\int_{B_{2/\rho}\setminus B_{1/\rho}} e^{-f}\eta^{\theta-2}v^{2-q}\le \epsilon\int_{\R^d} e^{-f}\eta^{\theta}v^{-q}+\frac{\rho^q}{\epsilon^{q-1} }\int_{B_{2/\rho}\setminus B_{1/\rho}} e^{-f}\eta^{\theta-q}.
%\eeq
From \eqref{sotto}, \eqref{young2volte} and by using \eqref{varphi}, we get
\begin{multline}\label{comp30}
(\frac{n}{2}+1-q)\int_{\R^d} e^{-f}v^{-q}|\nabla_{ g} v|_{ g}^2\varphi^\theta+\frac{2}{n-2}\int_{\R^d} e^{-f} v^{-q}\varphi^\theta\le  \theta\epsilon\int_{\R^d} e^{-f}\varphi^\theta v^{-q}|\nabla_{ g} v|_{ g}^2+ \theta \epsilon\int_{\R^d} e^{-f} v^{-q} \varphi^{\theta} \\
+\frac{\theta}{\epsilon^{q-1} R^q}\int_{B_{2R}^g\setminus B_R^g}  e^{-f}\varphi^{\theta-q}+\frac{\theta}{\epsilon^{q-1}r^q}\int_{B_{2r}^g\setminus B_{r}^g}  e^{-f}\varphi^{\theta-q}.
\end{multline}
Since $\theta > q$ then \eqref{VB_R} gives
$$
\frac{1}{r^q} \int_{B_{2r}^g\setminus B_{r}^g}  e^{-f}\varphi^{\theta-q} \leq \frac{1}{r^q}  \int_{B_{2r}^g\setminus B_{r}^g}  e^{-f} \leq C r^{n-q} \,.
$$
Since $\frac{n}{2}+1-q >0$ then $n-q>q-2\geq 0$ and this implies that 
$$
\lim_{r\to0} \frac{1}{r^q} \int_{B_{2r}^g\setminus B_{r}^g}  e^{-f}\varphi^{\theta-q} = 0.
$$
Hence, by taking the limit as $r\to 0$ in \eqref{comp30}, we obtain
\begin{equation}\label{comp30bis}
(\frac{n}{2}+1-q - \theta \epsilon)\int_{\R^d} e^{-f}v^{-q}|\nabla_{ g} v|_{ g}^2\varphi^\theta+(\frac{2}{n-2}- \theta \epsilon) \int_{\R^d} e^{-f} v^{-q}\varphi^\theta\le 
\frac{\theta}{\epsilon^{q-1} R^q}\int_{B_{2R}^g\setminus B_R^g}  e^{-f}\varphi^{\theta-q} \,.
\end{equation}
Since $\theta >q$, $(\frac{n}{2}+1-q)>0$ and by taking $\epsilon>0$ small enough, we obtain
$$
\int_{B_{2R}^g} e^{-f} v^{-q}|\nabla_{ g} v|_{ g}^2+\int_{B_{2R}^g}e^{-f}v^{-q}\le CR^{n-q},
$$
where $C=C(n,q,\epsilon)$. This concludes the proof for $2 \leq q<n/2+1$.

{\bf Step 2: proof of (ii)}. The case $q=0$ is exactly \eqref{VB_R}. If $2 \leq q<n/2+1$ then \eqref{(ii)} follows from \eqref{(i)}. 

Now we notice that if $0< q < 2$ then the assertion follows by a direct application of H\"older inequality and part (i) of the assertion. Indeed, let $0< s \leq 2$. We apply H\"older inequality with exponents $(\frac{q}{s},\frac{q}{q-s})$, with $2< q<n/2+1$, and get
$$
\int_{B_{2R}^g}e^{-f}v^{-s}\le\bigg(\int_{B_{2R}^g} e^{-f}v^{-q}\bigg)^{s/q}\bigg(\int_{B_{2R}^g}e^{-f}\bigg)^{(q-s)/q}\le C(R^{n-q})^{s/q}R^{n-sn/q}=CR^{n-s},
$$
where in the last inequality we used \eqref{(i)} for $2< q<n/2+1$ and \eqref{VB_R}. This implies that \eqref{(ii)} holds for $0\le q<n/2+1$. 

Hence, it remains to prove \eqref{(ii)} for $q=n/2+1$. From Lemma \ref{Lemma1} and by using \eqref{varphi} we have 
$$
\frac{2}{n-2}\int_{\R^d} e^{-f} v^{-\frac{n}{2}-1}\varphi^\theta \le\frac{\theta}{R} \int_{B_{2R}^g\setminus B_R^g}e^{-f} v^{-\frac{n}{2}}|\nabla_{ g} v|_{{g}}+\frac{\theta}{r}\int_{B_{2r}^g\setminus B_{r}^g}e^{-f}  v^{-\frac{n}{2}}|\nabla_{ g} v|_{{g}} \,,
$$
and from the definition of $\eta^\theta$ we obtain that
\beq\label{comp9}
\int_{B_{2R}^g \setminus B_{r}^g} e^{-f} v^{-\frac{n}{2}-1} \le\frac{C}{R} \int_{B_{2R}^g\setminus B_R^g}e^{-f} v^{-\frac{n}{2}}|\nabla_{ g} v|_{{g}} +\frac{C}{r} \int_{B_{2r}^g\setminus B_{r}^g}e^{-f} v^{-\frac{n}{2}}|\nabla_{ g} v|_{{g}}.
\eeq

\noindent Now we estimate the right-hand side of the last inequality. By writing $v^{-\frac{n}{2}}= v^{-\frac{n}{4}-\frac12+\frac{\epsilon}{2}}v^{-\frac{n}{4}+\frac12-\frac{\epsilon}{2}}$ and using H\"older inequality we find
$$
\int_{B_{2R}^g\setminus B_R^g} e^{-f} v^{-\frac{n}{2}}|\nabla_{ g} v|_{{g}}\le \bigg(\int_{B_{2R}^g\setminus B_R^g}e^{-f}v^{\epsilon-\frac{n}{2}-1}|\nabla_{ g} v|_{{g}}^2\bigg)^{1/2}\bigg(\int_{B_{2R}^g\setminus B_R^g}e^{-f}v^{-\frac{n}{2}+1-\epsilon}\bigg)^{1/2}.
$$
By choosing $0<\epsilon<\frac{n-2}{2}$ on the r.h.s.  we can apply \eqref{(ii)} for $2< q<\frac{n}{2}+1$ on the first term, and for $0<q<\frac{n}{2}+1$ on the second term to get
\beq\label{comp10}
\int_{B_{2R}^g\setminus B_R^g}e^{-f} v^{-\frac{n}{2}}|\nabla_{ g} v|_{{g}}\le CR^{(n+\epsilon-\frac{n}{2}-1)/2}R^{(n-\epsilon-\frac{n}{2}+1)/2}=CR^\frac{n}{2}.
\eeq
Analogously, the same argument applies to the second term on the r.h.s. of \eqref{comp9}, but with $r$ in place of $R$, and we find
\beq\label{comp10bis}
\int_{B_{2r}^g\setminus B_r^g}e^{-f} v^{-\frac{n}{2}}|\nabla_{ g} v|_{{g}}\le Cr^\frac{n}{2}.
\eeq
From \eqref{comp9}, \eqref{comp10} and \eqref{comp10bis} we obtain
\begin{equation*}
\int_{B_{2R}^g \setminus B_{r}^g} e^{-f} v^{-\frac{n}{2}-1} \le C (R^{\frac{n}{2}-1} + r^{\frac{n}{2}-1}) \,,
\end{equation*}
and by letting $r \to 0$ and using $n>2$ we find
\begin{equation*}
\int_{B_{2R}^g } e^{-f} v^{-\frac{n}{2}-1} \le C R^{\frac{n}{2}-1} \,,
\end{equation*}
and the proof is complete.
\end{proof}

We are ready to give the proof of Theorem \ref{thm_main}. It is clear that, after writing 
\begin{equation*} 
v=u^{-\frac{2}{n-2}}  \,,
\end{equation*}
Theorem \ref{thm_main} is a straightforward consequence of the following theorem.

\begin{theorem}\label{MainTh}
Let $d\ge 2$ and let $v$ be a positive solution to \eqref{eq2} satisfying \eqref{v_Dg}. If $\alpha^2\le \frac{d-2}{n-2}$ and $\frac{5}{2}<n\le 5$,  then
\begin{equation*}
v(x)= \begin{cases} 
c_1+c_2|x|^{2\alpha} & \textmd{ if } \alpha \neq 1 \,,\\
c_1+c_2|x-x_0|^{2\alpha} & \textmd{ if } \alpha = 1 \,,
\end{cases}
\end{equation*}
for any $x\in \R^d$, for some $c_1,c_2\in\R$ and $x_0 \in \R^d$. 
\end{theorem}

\begin{proof} 
From Proposition \ref{mainineq} we have that
\begin{equation} \label{prop41bis}
\frac{1}{\theta^2}(t-1/2)^2 \int_{\R^d} e^{-f}P^{t-2}v^{2-n}|\nabla_{{g}} P|_{{g}}^2\varphi^\theta\le\frac{1}{R^2}\int_{B_{2R}^g\setminus B_R^g} e^{-f} \varphi^{\theta-2}v^{2-n}P^t+\frac{1}{r^2}\int_{B_{2r}^g\setminus B_{r}^g} e^{-f}\varphi^{\theta-2}v^{2-n}P^t \,.
\end{equation}
From now on we consider $t=\frac{1}{2}+\delta$, where $\delta \in (0,\frac{1}{2})$ is small and to be chosen later.

We shall first show that the last term on the RHS in \eqref{prop41bis} vanishes for $r\to0$, and then we will show that also the limit as $R\to \infty$ vanishes. This will be enough to prove that  
%Observe that if we prove the following identities:
%\beq\label{I}
%\frac{1}{\theta^2}(t-1/2)^2 \int_{\R^d} e^{-f}P^{t-2}v^{2-n}|\nabla_{{g}} P|_{{g}}^2\varphi^\theta\le O(1) +\frac{1}{r^2}\int_{B_{2r}\setminus B_{r}} e^{-f}\varphi^{\theta-2}v^{2-n}P^t,
%\eeq
%for $R\to\infty$, and 
%\beq\label{II}
%\int_{B_{2r}\setminus B_{r}} e^{-f}\varphi^{\theta-2}v^{2-n}P^t=O(r^2)
%\eeq
%for $r\to 0$, then we can conclude that
%$$
%\int_{\R^d} e^{-f}P^{t-2}v^{2-n}|\nabla_{{g}} P|_{{g}}^2\varphi^\theta=0
%$$
%namely that 
$$
\nabla_{ g}P = 0 \quad \textmd{ in } \R^d \,. 
$$
Then \eqref{eqfond} implies that $k[v]=0$ and from Lemma \ref{lemma_k} we conclude.

\medskip

\textbf{Step 1: $r\to 0$.} In this step we show that  
\begin{equation} \label{r_to_0}
\int_{B_{2r}^g\setminus B_{r}^g} e^{-f}\varphi^{\theta-2}v^{2-n}P^t \leq C r^n \quad \textmd{ as } r \to 0 \,. 
\end{equation}
Indeed, since $v$ is locally bounded and $1\le (\frac{n-2}{2}vP)^{1-t}$, we have that 
$$
\int_{B_{2r}^g\setminus B_{r}^g} e^{-f}\varphi^{\theta-2}v^{2-n}P^t\le C\int_{B_{2r}^g\setminus B_{r}^g} e^{-f}+ C\int_{B_{2r}^g\setminus B_{r}^g} e^{-f}|\nabla_{ g} v|^2_{ g} \leq C r^n+C\int_{B_{2r}^g} e^{-f}|\nabla_{ g} v|^2_{ g}\,,
$$
and from Proposition \ref{DEL} and since $\alpha \leq 1$ we conclude.

\medskip

\noindent \textbf{Step 2: $R\to \infty$.} 
We shall split the proof in two steps, according to the case $\frac{5}{2}<n<5$ and $n=5$.

\textbf{- Step 2.1: $\frac{5}{2}<n<5$.} 
We assume that $\frac{5}{2}<n<5$ and we prove  
\beq\label{I}
\int_{B_{2R}^g\setminus B_R^g} e^{-f}\varphi^{\theta-2}v^{2-n}P^t \leq CR^{\frac{3}{2}-\delta}
\eeq
where $C$ does not depend on $R$.

Since $1-t>0$ and $\frac{n-2}{2}vP\ge 1$ we have $ 1\le (\frac{n-2}{2}vP)^{1-t}$, which implies
$$
\int_{B_{2R}^g\setminus B_R^g} e^{-f}\varphi^{\theta-2}v^{2-n}P^t  \le C\int_{B_{2R}^g\setminus B_R^g}e^{-f}\varphi^{\theta-2}v^{2-n-t}(vP)\le  CR^{2-t},
$$
where in the last inequality we applied Lemma \ref{Lemma2} that holds for
$$
2\le n+t-2<\frac{n}{2}+1,
$$
namely for
$$
4-t\le n <6-2t.
$$
Instead, for $2+t< n<4-t$ and $\delta$ small enough, by writing $v^{2-n}= v^{2-n+t} v^{-t}$ and applying H\"older inequality with exponents $\bigg(\frac{1}{t},\frac{1}{1-t}\bigg)$, we obtain
\beq\label{comp13}
\int_{B_{2R}^g\setminus B_R^g} e^{-f}\varphi^{\theta-2}v^{2-n}P^t \le\bigg(\int_{B_{2R}^g}e^{-f}v^{-1}P\bigg)^{t}\bigg(\int_{B_{2R}^g} e^{-f}v^{\frac{2-n+t}{1-t}}\bigg)^{1-t} \leq C R^{2-t}
\eeq
since   
\beq\label{comp12}
2\le \frac{n}{2}+1 \quad \text{ and } \quad 0< \frac{n-2-t}{1-t} \le \frac{n}{2}+1,
\eeq
and we can apply Lemma \ref{Lemma2} on the first and second terms.  We notice that \eqref{comp12} holds for $2+t<n<4-t$ by choosing $\delta$ small enough. Hence \eqref{I} is proved for $\frac{5}{2}<n<5$. 

\textbf{Step 2.2: $n=5$.}  Now we consider the case $n=5$ and we prove
\begin{equation} \label{In5}
\int_{B_{2R}^g\setminus B_R^g}e^{-f}\varphi^{\theta-2}v^{-3}P^t\le\epsilon_0^{-1}CR^{\frac{3}{2}-\delta}+C\epsilon_0 R^{2}\int_{B_{2R}^g\setminus B_R^g}e^{-f}\varphi^{\theta}v^{2-n}P^{t-2}|\nabla_{g}P|_{g}^2 \,,
\end{equation}
where $\delta>0$ is small enoug, $\epsilon_0 >0$ and $C$ does not depend on $\epsilon_0$ and $R$.
%
%We first notice that in this case Proposition \ref{mainineq} gives
%\beq\label{comp32}
%\frac{1}{\theta^2}(t-1/2)^2 \int_{\R^d} e^{-f}P^{t-2}v^{-3}|\nabla_{{g}} P|_{{g}}^2\varphi^\theta\le\frac{1}{R^2}\int_{B_{2R}^g\setminus B_R^g} e^{-f}\varphi^{\theta-2}v^{-3}P^t+\frac{1}{r^2}\int_{B_{2r}^g\setminus B_{r}^g} e^{-f}\varphi^{\theta-2}v^{-3}P^t.
%\eeq

We notice that, from the definition of $P$ \eqref{P-fun}, it follows 
\beq\label{pineq0}
1\le\frac{3}{2} vP \,,
\eeq
\begin{equation}\label{pineq}
|\nabla_{ g} v|_{{g}}^2\le \frac{5}{2}vP \,,
\end{equation}
and
\begin{equation}\label{pineq2}
vP \leq \frac{5}{2\epsilon}\bigg(\epsilon|\nabla_{ g} v|_{{g}}^2+\frac{2}{3}\bigg) \,,
\end{equation}
with $\epsilon \in (0,1]$.
From \eqref{pineq0} and recalling that $t= \frac12 + \delta$, we have
\begin{equation} \label{comp23}
\int_{B_{2R}^g\setminus B_R^g}e^{-f}\varphi^{\theta-2}v^{-3}P^t\le \bigg(\frac{3}{2}\bigg)^{\frac12 + \delta}\int_{B_{2R}^g\setminus B_R^g}e^{-f}\varphi^{\theta-2}v^{-\frac{5}{2}+\delta}P^{1+2\delta}
\end{equation}
and notice that an application of \eqref{pineq2} with $\epsilon=1$ gives
\begin{equation}\label{comp23bis}
\int_{B_{2R}^g\setminus B_R^g}e^{-f}\varphi^{\theta-2}v^{-\frac{5}{2}+\delta}P^{1+2\delta}
\le \frac{5}{2}  \int_{B_{2R}^g\setminus B_R^g}e^{-f}\varphi^{\theta-2}v^{-\frac{7}{2}+\delta}P^{2\delta}\bigg(|\nabla_{ g} v|_{{g}}^2+\frac{2}{3}\bigg) \,.
\end{equation}
We estimate the last term in \eqref{comp23bis} by using Corollary \ref{Corollary} with $q=\frac{7}{2} -\delta$, $m= 2\delta$ and $\psi=\varphi^{\theta-2}$, and we get  
\begin{multline}\label{comp14}
c(\theta,\delta)\int_{B_{2R}^g\setminus B_R^g}e^{-f}\varphi^{\theta-2}v^{-\frac{5}{2}+\delta}P^{1+2\delta} \\ 
\le \int_{B_{2R}^g\setminus B_R^g}e^{-f}\varphi^{\theta-3}v^{-\frac{5}{2}+\delta}P^{2\delta} | \nabla_{ g} v\cdot_{{g}}\nabla_{{g}}\varphi | 
+\int_{B_{2R}^g\setminus B_R^g}e^{-f}\varphi^{\theta-2}v^{-\frac{5}{2}+\delta}P^{-1+2\delta} | \nabla_{ g} v \cdot_{ g}\nabla_{ g}P | \\ = : J_1+J_2
\end{multline}
Now we focus on the first term on the RHS: $J_1$. An application of Young's inequality, with $\epsilon_1>0,q>1$ and exponents $\bigg(q,\frac{q}{q-1}\bigg)$ gives
$$
\varphi^{-1} P^{2\delta}   | \nabla_{ g} v\cdot_{{g}}\nabla_{{g}}\varphi |   \leq P^{1+ 2\delta}  \varphi^{-1} P^{-1}| \nabla_{ g} v|_g |\nabla_{{g}}\varphi |_g  \leq   C P^{1+ 2\delta}  \left\{ \frac{1}{\epsilon_1^{q-1}} \varphi^{-q} P^{-q} | \nabla_{ g} v|_g^q |\nabla_{{g}}\varphi |_g^q  + \epsilon_1   \right\} \,,
$$
and then $J_1$ in \eqref{comp14} can be bounded as follows
\begin{multline}\label{comp161}
J_1\le  \frac{C}{\epsilon_1^{q-1}} \int_{B_{2R}^g\setminus B_R^g}e^{-f}\varphi^{\theta-2-q}v^{-\frac{5}{2}+\delta}P^{1-q+2\delta}|\nabla_{ g} v|_{{g}}^q|\nabla_{{g}}\varphi|_{{g}}^q \\ +C \epsilon_1 \int_{B_{2R}^g\setminus B_R^g}e^{-f}\varphi^{\theta-2}v^{-\frac{5}{2}+\delta}P^{1+2\delta}=J_{11}+J_{12} \,.
\end{multline}
From \eqref{pineq} and \eqref{varphi} we obtain 
$$
J_{11}\le \frac{C}{\epsilon_1^{q-1}} R^{-q}\int_{B_{2R}^g\setminus B_R^g}e^{-f}\varphi^{\theta-2-q}v^{-\frac{5-q}{2}+\delta}P^{1-\frac{q}{2}+2\delta} \,.
$$
By choosing $q=2(1+2\delta)$ and applying Lemma \ref{Lemma2}, we have
$$
J_{11}\le \frac{C}{\epsilon_1^{1+4\delta}} R^{-2(1+2\delta)}\int_{B_{2R}^g\setminus B_R^g}e^{-f}\varphi^{\theta-4-4\delta}v^{-\frac{3}{2} + 3\delta} \leq  \frac{C}{\epsilon_1^{1+4\delta}}R^{ \frac32 - \delta } \,.
$$ 
Thus we get 
\beq\label{comp22}
J_1\le \frac{C}{\epsilon_1^{1+4\delta}}  R^{ \frac32 - \delta } +  \epsilon_1 C\int_{B_{2R}^g\setminus B_R^g}e^{-f}\varphi^{\theta-2}v^{-\frac{5}{2}+\delta}P^{1+2\delta}.
\eeq

Now we consider $J_2$ in \eqref{comp14}. From Young's inequality, we get
\begin{multline}\label{comp20}
J_2\le  \frac{1}{4\epsilon_0 R^2} \int_{B_{2R}^g\setminus B_R^g}e^{-f}\varphi^{\theta-4}v^{-2+2\delta}P^{-\frac{1}{2}+3\delta}|\nabla_{ g} v|^2_{{g}}
\\ +\epsilon_0 R^{2}\int_{B_{2R}^g\setminus B_R^g}e^{-f}\varphi^{\theta}v^{-3}P^{-\frac{3}{2} + \delta}|\nabla_{g}P|_{g}^2
=J_{21}+J_{22} \,.
\end{multline}
Focusing on the first term, by applying \eqref{pineq} we have
$$
J_{21} \leq \frac{1}{4\epsilon_0 R^2} \int_{B_{2R}^g\setminus B_R^g}e^{-f}\varphi^{\theta-4}v^{-1+2\delta}P^{\frac{1}{2}+3\delta}
$$
and by using H\"older inequality we obtain
\beq\label{comp33}
J_{21} \leq \frac{1}{4\epsilon_0 R^2} \bigg(\int_{B_{2R}^g\setminus B_R^g}e^{-f}\varphi^{\theta-4}v^{(-1+3\delta)\frac{3}{2}}P^{(\frac{1}{2}+3\delta)\frac{3}{2}}\bigg)^\frac{
 2}{3}\bigg(\int_{B_{2R}^g\setminus B_R^g}e^{-f}\varphi^{\theta-4}v^{-3\delta}\bigg)^\frac{ 1}{3} \,.
\eeq
Now from \eqref{pineq0} we have that
\beq\label{comp34}
1\le (\frac{n-2}{2}vP)^{1-(\frac{1}{2}+3\delta)\frac{3}{2}},
\eeq
then from \eqref{comp33} we get that
\beq\label{comp34}
J_{21} \leq \frac{C}{\epsilon_0 R^2} \bigg( \int_{B_{2R}^g\setminus B_R^g}e^{-f}\varphi^{\theta-4}v^{-\frac{9}{4}}vP\bigg)^\frac{
 2}{3}\bigg(\int_{B_{2R}^g\setminus B_R^g}e^{-f}\varphi^{\theta-4}v^{-3\delta}\bigg)^\frac{ 1}{3}.
\eeq
Since $n=5$, we can apply Lemma \ref{Lemma2} to \eqref{comp34} and obtain
\beq\label{comp35}
J_{21}\le \frac{C}{\epsilon_0}  R^{\frac{3}{2}-\delta}.
\eeq
Then from \eqref{comp20},\eqref{comp35} we obtain that
\beq\label{comp21}
J_2\le \frac{C}{\epsilon_0}  R^{\frac{3}{2}-\delta} + \epsilon_0 R^{2}\int_{B_{2R}^g\setminus B_R^g}e^{-f}\varphi^{\theta}v^{-3}P^{-\frac32 + \delta}|\nabla_{g}P|_{g}^2.
\eeq
Putting together \eqref{comp14}, \eqref{comp22} and \eqref{comp21} we have 
\begin{multline*}%\label{comp31}
\int_{B_{2R}^g\setminus B_R^g}e^{-f}\varphi^{\theta-2}v^{-\frac{5}{2}+\delta}P^{1+2\delta}\le C \left(\epsilon_1^{-(1+4\delta)} + \epsilon_0^{-1} \right)  R^{ \frac32 - \delta } \\
+  \epsilon_1 C\int_{B_{2R}^g\setminus B_R^g}e^{-f}\varphi^{\theta-2}v^{-\frac{5}{2}+\delta}P^{1+2\delta} + \epsilon_0 R^{2}\int_{B_{2R}^g\setminus B_R^g}e^{-f}\varphi^{\theta}v^{-3}P^{-\frac32 + \delta}|\nabla_{g}P|_{g}^2,
\end{multline*}
namely, for $\epsilon_1$ small enough 
\beq\label{comp36}
\int_{B_{2R}^g\setminus B_R^g}e^{-f}\varphi^{\theta-2}v^{-\frac{5}{2}+\delta}P^{1+2\delta} \le \epsilon_0^{-1}CR^{\frac{3}{2}-\delta}+C\epsilon_0 R^{2}\int_{B_{2R}^g\setminus B_R^g}e^{-f}\varphi^{\theta}v^{-3}P^{t-2}|\nabla_{g}P|_{g}^2,
\eeq
where $C=C(\theta, \delta,\epsilon_1)$. Hence, \eqref{comp23} and \eqref{comp36} yield \eqref{In5}.
%\begin{equation} \label{comp36bis}
%\int_{B_{2R}^g\setminus B_R^g}e^{-f}\varphi^{\theta-2}v^{-3}P^t\le\epsilon_0^{-1}CR^{\frac{3}{2}-\delta}+C\epsilon_0 R^{2}\int_{B_{2R}^g\setminus B_R^g}e^{-f}\varphi^{\theta}v^{2-n}P^{t-2}|\nabla_{g}P|_{g}^2 \,,
%\end{equation}

\textbf{Step 3: conclusion.} From \eqref{prop41bis}, \eqref{r_to_0}, \eqref{I} and \eqref{In5} (possibly by choosing $\epsilon_0$ small enough when $n=5$), we have
\begin{equation*}
\frac{\delta^2}{\theta^2} \int_{\R^d} e^{-f}P^{t-2}v^{2-n}|\nabla_{{g}} P|_{{g}}^2\varphi^\theta\le CR^{-\frac{1}{2}-\delta}  + o(1) \,,
\end{equation*}
as $r\to 0$. By letting $r\to0$ and then $R\to \infty$ we obtain that 
$$
\int_{\R^d} e^{-f}P^{t-2}v^{2-n}|\nabla_{{g}} P|_{{g}}^2  =0 \,.
$$
Hence $P$ is constant in $\R^d$. Hence Lemma \ref{lemma_P} implies that $k[v]=0$ and the conclusion follows from Lemma \ref{lemma_k}.
\end{proof}

\begin{proof}[Proof of Theorem \ref{thm_main}]
Theorem \ref{thm_main} is a straightforward consequence of Theorem \ref{MainTh}, recalling that $v=u^{-\frac{2}{n-2}}$ (see \eqref{v_def}).
\end{proof}

\section{Classification in cones} \label{section5}
In this section, we prove Theorem \ref{thm_main2}. Since the main arguments are strictly close to the ones used for the proof of Theorem \ref{thm_main}, we report here only the points where the presence of the cone plays some role. 

Let $A\subset \mathbb{S}^{d-1}$ be a smooth open domain in the $(d-1)$-dimensional sphere and consider the cone $\Sigma$ defined by
$$
\Sigma = \{x \in \R^d \,:\ x=t \xi \ \textmd{for some } t \in (0,+\infty) \textmd{ and } \xi \in A \} \,.
$$ 
We recall that in Theorem \ref{thm_main2} we assume that $\Sigma$ is convex and that we are considering positive solutions to 
\begin{equation}\label{eq_cone}
\begin{cases}
\divergence \left( |x|^{-2a} Du \right) + |x|^{-bp} u^{p-1} = 0 & \textmd{ in } \Sigma \\
u_\nu = 0  & \textmd{ on } \partial \Sigma \,,
\end{cases}
\end{equation}
with $u>0$ and $u \in D^{a,b}_{loc} (\overline\Sigma) $, where
\begin{equation*}
\mathcal{D}^{a,b}_{loc} (\overline\Sigma) : = \left\{ u \in L^p_{loc}(\overline \Sigma,|x|^{-b}dx) \textmd{ such that  } |x|^{-a} |Du| \in L^2_{loc}(\overline\Sigma,dx) \right\} \,.
\end{equation*} 
The weak formulation of \eqref{eq_cone} is the following
$$
\int_{\bar{\Sigma}} \nabla_{g} u \cdot_{g}\nabla_{g} \varphi \, e^{-f} dV_g = \int_{\bar{\Sigma}} u^{\frac{n+2}{n-2}} \varphi  \, e^{-f} dV_g \,=0,
$$
for any $\varphi \in C_c^\infty (\bar{\Sigma})$.

\begin{proof}[Proof of Theorem \ref{thm_main2} -- main remarks]  $ $ 

\noindent -- Formula \eqref{intbyparts1}. To obtain \eqref{intbyparts1} we integrated by parts. Due to the presence of the cone, in this case we obtain an additional boundary term which can be bounded thanks to the convexity of the cone. Indeed, we denote by $\text{II}(\nabla^Tv,\nabla^T v)$ the second fundamental form of $\partial \Sigma$ and we notice that
$$
 \nabla_{{g}} P \cdot_g \nu = \bigg(\frac{n}{2v}\nabla_g|\nabla_gv|^2-\frac{1}{v}P\nabla_gv\bigg)  \cdot_g \nu = -\frac{n}{2v}\text{II}(\nabla^Tv,\nabla^T v) \leq 0 \,,
$$
where we used the Neumann boundary condition and the convexity of the cone. Hence, since
\begin{multline*}
\left(t-\frac12\right) \int_{\Sigma}  e^{-f}P^{t-2}v^{2-n}|\nabla_{{g}} P|_{{g}}^2\varphi^\theta \\
\le - \theta\int_{\Sigma} \varphi^{\theta-1} e^{-f}v^{2-n}P^{t-1}\nabla_{{g}} P \cdot_g \nabla_{{g}}\varphi + \int_{\partial \Sigma} \varphi^{\theta} e^{-f}v^{2-n}P^{t-1}\nabla_{{g}} P \cdot_g \nu \,,\\
\end{multline*}
we obtain 
\begin{equation*}
\left(t-\frac12\right) \int_{\Sigma}  e^{-f}P^{t-2}v^{2-n}|\nabla_{{g}} P|_{{g}}^2\varphi^\theta \leq - \theta\int_{\Sigma} \varphi^{\theta-1} e^{-f}v^{2-n}P^{t-1}\nabla_{{g}} P \cdot_g \nabla_{{g}}\varphi 
\end{equation*}
which implies 
\begin{equation*}
\left(t-\frac12\right) \int_{\Sigma}  e^{-f}P^{t-2}v^{2-n}|\nabla_{{g}} P|_{{g}}^2\varphi^\theta \le \theta\bigg(\int_{\Sigma} e^{-f} 	\varphi^\theta v^{2-n}P^{t-2}|\nabla_{{g}} P|_{{g}}^2\bigg)^\frac{1}{2}\bigg(\int_{\Sigma} e^{-f}\varphi^{\theta-2}v^{2-n}P^t|\nabla_{{g}} \varphi|_{{g}}^2\bigg)^\frac{1}{2}, 
\end{equation*}
which is the analogue of \eqref{intbyparts1} in the case of a convex cone.

\noindent  -- Proof of Lemma \ref{Lemma1}. In the first step of the proof of Lemma \ref{Lemma1} we integrate by parts.  The boundary terms that we have in the case of the cone are still zero thanks to the boundary conditions.

\noindent -- Lemma \ref{lemma_k}. The point $x_0$ appearing in \eqref{v_radial} needs another condition, in particular, $x_0$ must be a vertex of the cone. This restriction follows from the Neumann boundary condition: since $\nabla_g v(x)$ points in the direction of $x_0 -x$ then the Neumann boundary condition implies that $x_0$ must be a vertex.

%  namely we know that the direction of $\nabla v(x)$ is given by $x-x_0$. 
%\begin{itemize}
%\item If $x_0\in \R^n\setminus\Gamma$ consider the ball centered in $x_0$ such that is tangent to the cone in a point $\bar x\in \Gamma$.  Then the gradient in $\bar x$ has the same direction as the outward normal in $\bar x$, namely $\partial_\nu v(\bar x)\ne 0$ with $\bar x\in \Gamma$, then $v$ isn't a solution.
%\item If $x_0\in\Gamma$ but it isn't a vertex, consider the ball centred in $x_0$ such that is tangent to the cone in another point $\tilde x\in \Gamma$.  From the construction of the cone, $\tilde{x}$ can't be a vertex.  Then the gradient in $\tilde x$ has the same direction as the outward normal in $\tilde x$,  namely $\partial_\nu v(\tilde x)\ne 0$ with $\tilde x\in \Gamma$, then $v$ isn't a solution.
%\item If $x_0$ is a vertex and we consider a generic point $x'$ of $\Gamma$, then $\nabla v(x')$ has the same direction of $x'-x_0$ that is orthogonal to the normal in $x'$ by construction of the cone, then $v$ is a solution. 
%\end{itemize}

\noindent -- Proposition \ref{DEL}. The proof of this proposition strongly uses the results in \cite[Propositions 8.1 and 8.2]{DEL}. A careful analysis of the proof in \cite{DEL} reveals that all the arguments can be repeated almost verbatim in the case of a convex cone. Indeed, Kelvin and Emden-Fowler transformations still apply in this case and, by setting 
$$
\varphi (s,\omega)=r^{a_c-a}v(r,\omega)\qquad s=-\log r, \omega=\frac{x}{r}\,,
$$
one has to deal with a problem for $\varphi$ and $\tilde \varphi$ (which is the Kelvin transform of $\varphi$) in a cylinder, which has the form 
\begin{equation*}
\begin{cases}
-\partial^2_s\varphi-\Delta_\omega\varphi+\Lambda\varphi=\varphi^{p-1} & \textmd{ in } \mathbb R \times A \\
\partial_\nu \varphi = 0  & \textmd{ on } \mathbb R \times \partial A \,,
\end{cases}
\end{equation*}
where $\Lambda=(a_c-a)^2$. While the presence of the Neumann boundary condition does not introduce any obstacle for repeating the argument in \cite[Propositions 8.1 and 8.2]{DEL}, the convexity of $\Sigma$ is a crucial condition. Indeed, since $\Sigma$ is convex then the first nontrivial Neumann eigenvalue of $\Delta_\omega$ in $A$ is greater or equal $d-1$ (see \cite[Theorem 4.3]{Escobar}), which is the threshold for succeeding in proving the results.  
\end{proof}

\appendix

\section{A differential identity in Riemannian Manifolds}

In this appendix, we prove a differential identity which is crucial for our approach. This identity was essentially proved in \cite[Formula (2.18)]{CCR} (see also \cite{CFP}), but it holds for more general weights and Riemannian manifolds. For this reason and in order to make the paper self-contained, we prefer to state and prove it in the following lemma.

\begin{lemma} \label{lemma_app}
Let $(M,g)$ be a Riemannian manifold of dimension $d$. Let $n \in \mathbb{R}$, $f \in C^2(M)$ and $v \in C^3(M)$, with $v>0$. Then we have
\begin{equation*}
e^f \divergence \bigg( e^{-f} v^{1-n}\bigg( \frac{n}{2}\nabla_{{g}}|\nabla_{{g}} v|^2_{{g}} - n \Cl v\nabla_{{g}} v + (n-1)P \nabla_{{g}} v\bigg)\bigg) \\
 =n v^{1-n}k[v] + v \bigg(\Cl v- P \bigg) \Cl(v^{1-n}),
\end{equation*}
with
\begin{equation*} 
\Cl v = e^f \divergence_{g}\bigg(e^{-f} \nabla_{{g}}v\bigg) 
\end{equation*}
and
\begin{equation*}
P=\frac{2}{n-2}\frac{1}{v}+\frac{n}{2}\frac{|\nabla_{{g}} v|^2_{g}}{v} \,,
\end{equation*}
and where we set
\begin{equation*} 
k[v]=\big{|}H_{{g}} v\big{|}_g^2 -\frac{(\Cl v)^2}{n} +Ric_{{g}}(\nabla_{{g}} v,\nabla_{{g}} v)+H_{{g}}f(\nabla_{{g}} v,\nabla_{{g}} v) \,.
\end{equation*}
\end{lemma}

\begin{proof}
To simplify the notation, in this proof, we omit the subscript $g$. We first notice that $\Cl v$ can be written as 
\begin{equation} \label{app0.1}
\Cl v = \Delta v - \nabla f \cdot \nabla v
\end{equation}
and then 
\begin{equation*}
\begin{split}
\frac{n}{2} e^f \divergence & \left( e^{-f}  v^{1-n}  \nabla |\nabla v|^2  \right)   = \frac{n}{2} v^{1-n} \Cl \left(  |\nabla v|^2   \right) + (1-n) \frac{n}{2} v^{-n} \nabla v \cdot \nabla |\nabla v|^2 \\
& = \frac{n}{2} v^{1-n} \Delta \left( |\nabla v|^2   \right) - \frac{n}{2} v^{1-n}  \nabla f \cdot \nabla |\nabla v|^2 + (1-n) \frac{n}{2} v^{-n} \nabla v \cdot \nabla |\nabla v|^2  \,,
\end{split}
\end{equation*}
and Bochner identity yields
\begin{multline} \label{app1.1}
\frac{n}{2} e^f \divergence \left( e^{-f}  v^{1-n}  \nabla |\nabla v|^2  \right)  \\ = n v^{1-n}  \left( |H v|^2 + \nabla v \cdot \nabla \Delta v  + Ric( \nabla v ,\nabla v)   \right) - \frac{n}{2} v^{1-n}  \nabla f \cdot \nabla |\nabla v|^2 + (1-n) \frac{n}{2} v^{-n} \nabla v \cdot \nabla |\nabla v|^2 \,.
\end{multline}
Now we use \eqref{app0.1} to write
\begin{equation*}
\begin{split}
&-n  e^f \divergence  \left( e^{-f} v^{1-n}  \Cl v\nabla v  \right)  = -n \left[ v^{1-n} (\Cl v)^2 + v^{1-n} \nabla v \cdot \nabla \Cl v + (1-n)v^{-n} \Cl v |\nabla v|^2    \right]  \\
 & = -n \left[ v^{1-n} (\Cl v)^2 + v^{1-n} \nabla v \cdot \nabla \Delta v - v^{1-n} \nabla v \cdot \nabla (\nabla f \cdot \nabla v)    + (1-n)v^{-n} \Cl v |\nabla v|^2   \right] \\
 & = -n \left[ v^{1-n} (\Cl v)^2 + v^{1-n} \nabla v \cdot \nabla \Delta v - v^{1-n} Hf (\nabla v,\nabla v) - v^{1-n} \nabla f \cdot \nabla \frac{|\nabla v|^2}{2}   + (1-n)v^{-n} \Cl v |\nabla v|^2   \right] \,,
\end{split}
\end{equation*}
and from \eqref{app1.1} we have 
\begin{multline*}
e^f \divergence \left\{ e^{-f} \left( \frac{n}{2} v^{1-n}  \nabla |\nabla v|^2  -nv^{1-n}  \Cl v\nabla v  \right) \right\}
= n v^{1-n}  \left[ |H v|^2 - (\Cl v)^2   + Ric( \nabla v ,\nabla v) + Hf (\nabla v,\nabla v)  \right] \\
+ (1-n) v^{-n} \left(  \frac{n}{2} \nabla|\nabla v|^2 \cdot \nabla v  -n \Cl v |\nabla v|^2   \right) \,,
\end{multline*}
i.e. 
\begin{multline} \label{app1.2}
e^f \divergence \left\{ e^{-f} \left( \frac{n}{2} v^{1-n}  \nabla |\nabla v|^2  -nv^{1-n}  \Cl v\nabla v  \right) \right\}
= n v^{1-n} k[v] + (1-n)  v^{1-n}   (\Cl v)^2  \\
+ (1-n) v^{-n} \left(  \frac{n}{2} \nabla|\nabla v|^2 \cdot \nabla v  -n \Cl v |\nabla v|^2   \right) \,.
\end{multline}
%and, since $ \frac{n}{2}\nabla  |\nabla v|^2 = \nabla(vP) = P \nabla v + v \nabla P$ and by the definition of $k[v]$, we obtain
%\begin{multline*}
%\Cl \left( \frac{n}{2} v^{1-n}  \nabla |\nabla v|^2  -nv^{1-n}  \Cl v\nabla v  \right) 
%= n v^{1-n} k[v] + (1-n)  v^{1-n}   (\Cl v)^2 + \\
%+ (1-n) v^{-n} \left( P |\nabla v |^2   -n \Cl v |\nabla v|^2   \right)  + (1-n) v^{1-n} \nabla P \cdot \nabla v \,.
%\end{multline*}
%
%
Now we consider
\begin{equation*}
\begin{split}
(n-1) e^f \divergence & (e^{-f} v^{1-n} P \nabla v) \\
&  = (n-1) v^{1-n} \left( P \Cl v  + \nabla P \cdot \nabla v \right) - (n-1)^2 v^{-n}  P |\nabla v|^2  \\
& =  (n-1) v^{1-n} \left( P \Cl v  - \frac{P}{v} |\nabla v|^2 + \frac{n}{2v} \nabla |\nabla v|^2 \cdot \nabla v \right) - (n-1)^2 v^{-n}  P |\nabla v|^2\,,
\end{split}
\end{equation*}
where in the last equality we used  
$$
\nabla P = - \frac{P}{v} \nabla v + \frac{n}{2v} \nabla |\nabla v|^2\,,
$$
and hence
\begin{multline} \label{app1.3}
(n-1) e^f \divergence (e^{-f} v^{1-n} P \nabla v) \\ =  (n-1) v^{1-n}  P \Cl v   +  \frac{n}{2}  (n-1) v^{-n} \nabla |\nabla v|^2 \cdot \nabla v   - n(n-1) v^{-n}  P |\nabla v|^2 \,,
\end{multline}
Hence from \eqref{app1.2} and \eqref{app1.3} we find
\begin{multline*}
e^f \divergence \left( e^{-f}  v^{1-n} \left[ \frac{n}{2} \nabla |\nabla v|^2  -n  \Cl v\nabla v + (n-1)  P \nabla v \right] \right) 
= n v^{1-n} k[v] + (1-n)  v^{1-n}   (\Cl v)^2  \\
+ n (n-1) v^{-n}  \Cl v |\nabla v|^2 + (n-1) v^{1-n} P   \Cl v  - n (n-1) v^{-n} P  |\nabla v|^2   \,,
 \end{multline*}
i.e.
\begin{multline*}
e^f \divergence \left( e^{-f} v^{1-n} \left[ \frac{n}{2} \nabla |\nabla v|^2  -n   \Cl v\nabla v + (n-1)   P \nabla v \right] \right) \\
= n v^{1-n} k[v]  + n(n-1) v^{-n}   |\nabla v|^2 (\Cl v - P) + (n-1) v^{1-n}  \Cl v (P- \Cl v)
\end{multline*}
and thus
\begin{multline*}
e^f \divergence \left( e^{-f} v^{1-n} \left[ \frac{n}{2} \nabla |\nabla v|^2  -n   \Cl v\nabla v + (n-1)   P \nabla v \right] \right) \\
= n v^{1-n} k[v]  + (n-1) v^{-n} (n  |\nabla v|^2 - v \Cl v)(\Cl v - P)   \,.
\end{multline*}
Since
$$
v \Cl (v^{1-n}) = (n-1) v^{-n} \left[ n  |\nabla v|^2 - v \Cl v \right]
$$
we conclude.
\end{proof}

\section*{Acknowledgments} 
\noindent The authors thank Alberto Farina for useful discussions and remarks.
The authors have been partially supported by the ``Gruppo Nazionale per l'Analisi Matematica, la Probabilit\`a e le loro Applicazioni'' (GNAMPA) of the ``Istituto Nazionale di Alta Matematica'' (INdAM, Italy) and by the Research Project of the Italian Ministry of University and Research (MUR) Prin 2022 
``Partial differential equations and related geometric-functional inequalities'', grant number 20229M52AS\_004.

\end{document}